\documentclass{amsart}
\numberwithin{equation}{section}

\newtheorem{tl}{Corollary}
\newtheorem{dl}{Theorem}[section]
\newtheorem{yl}{Lemma}
\newtheorem{dy}{Definition}
\newtheorem{lz}{Example}
\newtheorem{xz}{Proposition}

\newtheorem{rmm}{Remark}

\theoremstyle{remark}

\def\qed{\hfill \rule{4pt}{7pt}}
\def\pf{\noindent {\it Proof.} }
\newcommand{\poq}[2]{(#1;q)_{#2}}
\newcommand{\poqp}[2]{(#1;p)_{#2}}
\newcommand{\poqt}[2]{(#1;|q|)_{#2}}
\newcommand{\poqpt}[2]{(#1;|p|)_{#2}}
\begin{document}
\title[Generalized Ismail's argument and $(f,g)$-expansion Formula]{Generalized
Ismail's argument and\\ $(f,g)$-expansion Formula}
\author{X. Ma}
\date{2.11.2006}
\address{Department of Mathematics, SuZhou University, SuZhou,P. R. China}

\email{xrma@public1.sz.js.cn}

\thanks{}
\date{8.25.2006}
\subjclass[2000]{Primary 05A10,05A19,33D15; Secondary
5A15,33C20,33D20}

\keywords{Lagrange inversion formula,  $(f,g)$-inversion,
$(f,g)$-difference operator, $(f,g)$-series, $q$-series, basic
hypergeometric series, divided difference operator, $q$-difference
operator, convergency problem, representation problem,
$(f,g)$-expansion formula, Ismail's argument, transformation
formula, summation formula.}

\begin{abstract}
As further development of earlier works on the $(f,g)$-inversion,
the present paper is devoted to
 the $(f,g)$-difference operator and the representation problem or an expansion formula of analytic functions.
  A recursive formula and
the Leibniz formula for the $(f,g)$-difference operator of the
product of two functions are established.
 The resulting expansion formula not only
unifies  the $q$-analogue of the Lagrange inversion formula of
Gessel and Stanton
 (thus, a $q$-expansion formula of Liu) for $q$-series but also systematizes the ``Ismail's
argument". In the meantime,  a rigorous analytic proof of the
$(1-xy,x-y)$-expansion formula with respect to
   geometric series, along with a proof of the previously unknown fact that it is equivalent
  to  a $q$-analogue of the
 Lagrange inversion formula due to Gessel
and Stanton, is presented. As applications, new proofs of several
well-known summation and transformation formulas are investigated.
\end{abstract}
\maketitle

\section{Introduction}
It is well known that the core of the classical Lagrange inversion
formula (cf.\cite[\S 7.32]{1001}) is to express  the coefficients
$a_n$ in the expansion of
\begin{eqnarray}
  F(x)=\sum_{n=0}^{\infty}a_n\left(\frac{x}{\phi(x)}\right)^n\label{lag}
\end{eqnarray}
by
$$
a_n=n!\frac{d^{n-1}}{dx^{n-1}}\left[\phi^{n}(x)\frac{dF(x)}{dx}\right]_{x=0}
$$
provided that $F(x)$ and $\phi(x)$ are analytic around $x=0$, $
\phi(0)\neq 0$, $\frac{d}{dx}$ denotes the usual derivative
operator.

 In the past years,
 various $q$-analogues (as generalizations) of the Lagrange inversion formula, as an
active field of research with an increasing number of applications
to $q$-series and the Rogers-Ramanujan identities, have been studied
by numerous authors (cf.\cite{andrew,6,111,kratt,18,22,11}). For a
good survey about
 results and open problems on this topic, we would like to refer the
reader to Stanton's paper \cite{11} and  only repeat here, for
reference purposes later, three noteworthy results.

One $q$-analogue  found by Carlitz in 1973
(cf.\cite[Eq.(1.11)]{carliz}), subsequently reproduced by Roman
(cf.\cite[p. 253, Eq. (8.4)]{3rr}) viz $q$-umbral calculus, is that
for any formal
 series $F(x)$, it holds that
\begin{eqnarray}
F(x)=\sum_{k=0}^{\infty}\frac{x^k}{(q,x;q)_k}\left[\mathcal{D}^k_{q,x}\{F(x)(x;q)_{k-1}\}\right]_{x=0},\label{carbe}
\end{eqnarray}
where  $\mathcal{D}_{q,x}$ denotes the $q$-difference operator.

Ten years later, Gessel and Stanton \cite{111} successfully,  with
insight
 that the essential character of the Lagrange inversion
formula is equivalent to finding the inverse of an infinite lower
triangular matrix $F=(B_{n,k})$  subject to $B_{n,k}=0$ unless $\
n\geq k$, $B_{n,n}\neq 0$, that is, another unique matrix
$G=(B_{n,k}^{-1})$ satisfying $$ \ {{\sum }_{n\geq i\geq
k}}B_{n,i}^{-1}B_{i,k}=\delta_{n,k},
$$
where $\delta$ denotes the usual Kronecker delta (as always,  such a
pair of $F$ and $G$ is called a {\sl matrix inversion}),  discovered
a few $q$-analogues of some matrix inversions. See \cite[Theorems
3.7 and 3.15]{111}. One result is  that for any formal power series
$F(x)$,
\begin{align} F(x)=\sum_{n\geq k\geq
0}^{\infty}a_k\frac{(Ap^kq^k;p)_{n-k}}{\poq{q}{n-k}}q^{-nk}x^n\label{gessel1}
\end{align}
if and only if
\begin{align}
a_n=\sum_{k=0}^n(-1)^{n-k}q^{\binom{n-k+1}{2}+nk}
\frac{(1-Ap^kq^k)(Aq^np^{n-1};p^{-1})_{n-k-1}}{\poq{q}{n-k}}F(q^k).\label{gessel2}
\end{align}

 In 2002, Liu \cite{liu} established by using of Carlitz's $q$-analogue   (\ref{carbe})
 and the special case $c\mapsto \infty$ of Rogers'  $\,_6\phi_5$ summation
formula \cite[II.20]{10}, as well as the $q$-difference operator,
the following $q$-expansion formula: given any formal series $F(x)$,
it holds that
\begin{eqnarray}
F(x)=\sum_{k=0}^{\infty}\frac{(1-aq^{2k})(aq/x;q)_kx^k}{(q;q)_k(x;q)_k}\left[\mathcal{D}^k_{q,x}\{F(x)(x;q)_{k-1}\}
\right]_{x=aq} \label{1111}.
\end{eqnarray}

As of today,  these expansion formulas have been proved to be very
important to the theory of  basic hypergeometric series (i.e.,
$q$-series), but they were established without considering
convergence and are therefore only valid within the ring of formal
power series. See the above references for more details. The
investigation of their rigorous analytic proofs, which is not only
unavoidable but also essential to basic hypergeometric series and
special function, is one of the purposes of the present paper.

 To achieve this goal, one important idea we will invoke, particular
 regarding elementary derivations of summation and transformation formulas of
basic hypergeometric series, is  Ismail's  analytic proof of
Ramanujan's $\,_1\varphi_1$ summation formula \cite{ismail0}. It was
often referred to as the ``Ismail's argument" (cf.\cite{schlosser1})
since it is Ismail who was apparently the first to apply analytic
continuation argument in the context of bilateral basic
hypergeometric series. Later, Askey and Ismail used this method  and
Rogers' $\,_6\phi_5$ summation formula \cite{ask} to evaluate
Bailey's $\,_6\varphi_6$ sum. In principle, the ``Ismail's argument"
can be summarized briefly as follows: if one wants to prove two
analytic functions $F(x)=G(x)$,  all that is necessary is to show
that they agree infinitely often near a point that is an interior
point of the set of analyticity. Unfortunately, application of this
idea to $q$-series has not been investigated systematically. Even
later, it was treated extensively by Gasper in \cite{188}, from
which one can see how the analytic continuation of a given function
affects existences of summation and transformation formulas of
$q$-series, but this method has not yet been written in their
remarkable book \cite{10}. One reason for this is that, adopting
Gasper as saying`` $\cdots$ the succeeding higher order derivatives
becomes more and more difficult to calculate for $|z|<1$, and so one
is forced to abandon this approach and to search for another way
$\cdots$" . See \cite{188} for more details.

The purposes of the present paper are: (1) to show that all above
$q$-analogues or expansion formulas can be unified by, from a purely
analytic viewpoint, the representation problem of $F(x)$ in terms of
a new kind of series, that is $(f,g)-$series (see Definition
\ref{deridef00} below),
 \begin{eqnarray}F(x)\sim\sum_{n=0}^{\infty}G(n)
f(x_n,b_n)
\frac{\prod_{i=0}^{n-1}g(b_i,x)}{\prod_{i=1}^{n}f(x_i,x)}\label{expandef00}
\end{eqnarray}
 where the coefficient $G(n)$ denotes the $n$-th order $(f,g)$-difference of
  $F(x)$ (see Definition \ref{deridef} below)
 expressed explicitly  by (\ref{expandef0}),
 $f(x,y)\in \mbox{\sl Ker}\mathcal{L}^{(g)}_{3}$ the set of functions over ${\Bbb C}\times{\Bbb C}$
 in two
variables
 $x,y$, ${\Bbb C}$ is the usual complex field,
such that for all $a,b,c,x\in \mathbb{{\Bbb C}}$
\begin{eqnarray}
f(x,a)g(b,c)+f(x,b)g(c,a)+f(x,c)g(a,b)=0,\label{idd}
\end{eqnarray}
where $g(x,y)$ is antisymmetric, i.e., $g(x,y)=-g(y,x)$. Later as we
will see, they  are special cases of the following expansion formula
\begin{eqnarray}\label{ggg} F(x)=\sum_{n=0}^{\infty}G(n)
f(x_n,b_n)
\frac{\prod_{i=0}^{n-1}g(b_i,x)}{\prod_{i=1}^{n}f(x_i,x)}\
\end{eqnarray}
under the assumption that the series in the r.h.s.\,of
(\ref{expandef00}) converges to $F(x)$. In a certain sense, it
provides a \emph{discrete} analogue of the Lagrange inversion
formula by replacing $(x/\phi(x))^n$ in (\ref{lag}) with
$\prod_{i=0}^{n-1}g(b_i,x)/\prod_{i=1}^{n}f(x_i,x)$. (2) to set up
an existence theorem of (\ref{ggg}) when $F(x)$, $\{b_i\}$, and
$\{x_i\}$ are subject to suitable convergence conditions. We believe
it can serve as a standard model for the ``Ismail's argument". We
call it the generalized ``Ismail's argument". (3) to apply the
generalized ``Ismail's argument" to basic hypergeometric series in
order to derive new or review known summation formulas.

 Our argument is based on a recent discovery that: as
our previously paper \cite{0020} shows, Identity (\ref{idd}) implies
the following matrix inversion, and vice versa. It is called the
($f,g$)-inversion formally in \cite{0021} by the author.

\begin{dl} \label{math4}
 Let $F=(B_{n,k})_{n,k\in {\Bbb Z}}$ and $G=(B^{-1}_{n,k})_{n,k\in {\Bbb Z}}$ be two matrices
with entries given by
\begin{eqnarray}
B_{n,k}&=&\frac{\prod_{i=k}^{n-1}f(x_i,b_k)}
{\prod_{i=k+1}^{n}g(b_i,b_k)}\label{news1115550}\qquad\mbox{and}\\
B^{-1}_{n,k}&=&
\frac{f(x_k,b_k)}{f(x_n,b_n)}\frac{\prod_{i=k+1}^{n}f(x_i,b_n)}
{\prod_{i=k}^{n-1}g(b_i,b_n)},\label{news1115551}\quad\mbox{respectively},
\end{eqnarray}
where ${\Bbb Z}$ denotes the set of integers, $\{x_i\}$ and
$\{b_i\}$ are arbitrary sequences such that none of the denominators
in the right hand sides of (\ref{news1115550}) and
(\ref{news1115551}) vanish. Then $F=(B_{n,k})_{n,k\in {\Bbb Z}}$ and
$G=(B^{-1}_{n,k})_{n,k\in {\Bbb Z}}$ is a matrix inversion if and
only if for all $a,b,c,x\in \mathbb{{\Bbb C}}$,  (\ref{idd}) holds.
\end{dl}
The reader may consult \cite{0020,0021} for  more detailed
expositions.

The present paper, as further study of the $(f,g)$-inversion, is
organized as follows. In Section 2, the $n$-th order
$(f,g)$-difference operator   is introduced and  some basic
propositions like the {\sl Leibniz formula} are examined. Section 3
 is devoted to the generalized ``Ismail's argument" and some $(f,g)$-expansion
 formulas. They are all special solutions of the representation problem
 of analytic functions in terms of $(f,g)$-series.
 A rigorous analytic proof of  the $(1-xy,x-y)$-expansion formula with respect to
  two  sequences of geometric series as well as a proof of the previously unknown fact that it is equivalent
  to  the $q$-analogue (\ref{gessel1})/(\ref{gessel2}) of the
 Lagrange inversion formula due to Gessel
and Stanton are given in Sections 4 and 7, respectively.   Their
applications to  basic hypergeometric series are investigated in
details in Section 5. In Section 6, a few concluding remarks on some
 problems related with the $(f,g)$-expansion formulas are given.

{\sl Notations and conventions.} Throughout this paper we adopt the
standard notation and terminology for basic hypergeometric series of
Gasper and Rahman's book \cite{10}, which is also a main reference
for any related results concerning basic hypergeometric series. For
instance, given a (fixed) complex number $q$ with $|q|<1$, a complex
$a$ and a natural number $n$, define the $q$-shifted factorials $(a;
q)_n$ and $(a;q)_\infty$  as
\begin{eqnarray*}
 (a;q)_\infty =
\prod_{k=0}^{\infty}(1-aq^k),\quad
(a;q)_n=(a;q)_\infty/(aq^n;q)_\infty
\end{eqnarray*}
with the following compact multi-parameter notation
\begin{eqnarray*}
  (a_1,a_2,\cdots,a_m;q)_n = (a_1;q)_n(a_2;q)_n\cdots
   (a_m;q)_n
    \end{eqnarray*}
and the $q$-binomial coefficient ${n\brack
k}_q=\poq{q}{n}/(\poq{q}{k}\poq{q}{n-k}). $   The {\em basic
hypergeometric series\/} with the base $q$ is defined by
\begin{equation*}{}_{r}\phi _{s}\left[\begin{matrix}a_{1},\dots ,a_{r}\\ b_{1},\dots ,b_{s}\end{matrix}
; q, z\right]=\sum _{n=0} ^{\infty }\frac{\poq {a_{1},\cdots
,a_{r}}{n}}{\poq
{q,b_{1},\cdots,b_{s}}{n}}\left[(-1)^nq^{n(n-1)/2}\right]^{1+s-r}z^{n},\end{equation*}
while the  \emph{bilateral basic hypergeometric series} is defined
as
\begin{equation*}{}_{r}\varphi_{s}\left[\begin{matrix}a_{1},\dots ,a_{r}\\ b_{1},\dots ,b_{s}\end{matrix}
; q, z\right]=\sum _{n=-\infty} ^{\infty }\frac{\poq {a_{1},\cdots
,a_{r}}{n}}{\poq
{b_{1},\cdots,b_{s}}{n}}\left[(-1)^nq^{n(n-1)/2}\right]^{s-r}z^{n}.\end{equation*}
In the setting of bilateral series, we employ the convention of
defining (a product form of Eq.(3.6.12) in \cite{10}) over  ${\Bbb
Z}$
\begin{eqnarray*}
\prod_{j=k}^mA_j= \left\{\begin{array}{ll} A_{k}A_{k+1}\cdots
A_{m},& m\geq k;\\
1,& m=k-1;\\
 (A_{m+1}A_{m+2}\cdots A_{k-1})^{-1},& m\leq k-2.
\end{array}
\right.
\end{eqnarray*}
\subsection*{Acknowledgements} The author would like to thank his
teacher L. Zhu and Prof. L.S.Wu  for deep and interesting
discussions and many inspiring ideas during the preparation of this
paper. He also wants to thank Z.G. Liu and R.M. Zhang for their
helpful discussions during the period of the Workshop  on $q$-series
held at ECNU in 2006.

\section{The $n$-th order $(f,g)$-difference operator}\setcounter{equation}{0}\setcounter{dl}{0}\setcounter{tl}{0}\setcounter{yl}{0}\setcounter{lz}{0}
Before we turn to our main results, we need the definition of
 $(f,g)$-difference operator which will play the same role in our
 expansion formulas as the divided difference in the
classical Newton's interpolation formula.  In what follows, we use
${\Bbb C}(x)$
 to denote the linear space of all functions  over $\mathbb{\Bbb C}$ of a variable
$x$.
\begin{dy}\label{ddd} Let $f(x,y)\in \mbox{\sl Ker}\mathcal{L}^{(g)}_{3}$, $\{x_i\}$ and $\{b_i\}$ be arbitrary sequences such that none of
the denominators in  (\ref{deridef}) vanish.
 Then the mapping $$ \,{\bf\mathbb{D}}^{(n)}_{(f,g)}\left[\begin{matrix}b_0,b_{1},\dots ,b_{n}\\
x_{1},\dots ,x_{n-1}\end{matrix}\right]\{\bullet\}: {\Bbb
C}(x)\longrightarrow {\Bbb C},$$ such that
\begin{eqnarray}
\,{\bf\mathbb{D}}^{(n)}_{(f,g)}\left[\begin{matrix}b_0,b_{1},\dots ,b_{n}\\
x_{1},\dots ,x_{n-1}\end{matrix}\right]\{F(x)\}=
\sum_{k=0}^{n}F(b_k)\frac{\prod_{i=1}^{n-1} f(x_i,b_k)}
 {\prod_{i=0,i\neq k}^{n}g(b_i,b_k)}\label{deridef}
\end{eqnarray}
 is said to be the $n$-th order $(f,g)$-difference
operator   with respect to $n+1$ pairwise distinct nodes
$b_0,b_1,b_2,$ $\cdots, b_n$ and $n-1$ parameters
$x_1,x_2,\cdots,x_{n-1}.$
\end{dy}
In what follows, we often abbreviate
$\,{\bf\mathbb{D}}^{(n)}_{(f,g)}\left[\begin{matrix}b_0,b_{1},\dots ,b_{n}\\
x_{1},\dots ,x_{n-1}\end{matrix}\right]\{\bullet\}$ by
$\,{\bf\mathbb{D}}^{(n)}_{(f,g)}$ if it is clear in the context.
 When $f=1$ or $n=1$, the notation
$\,{\bf\mathbb{D}}^{(n)}_{(g)}[b_0,b_1,b_2,\cdots,b_n]$ is employed
in which $\{x_i\}$ are bypassed, since in either of these two cases,
the sum does certainly not depend on $\{x_i\}$. In
particular, for $n=0$, by the convention,  $$\,{\bf\mathbb{D}}^{(n)}_{(f,g)}\left[\begin{matrix}b_0\\
x_{-1}\end{matrix}\right]\{F(x)\}=\frac{F(b_0)}{f(x_0,b_0)}.$$

 The following example displays
that $\,{\bf\mathbb{D}}^{(n)}_{(f,g)}$ is a generalization of  the
divided difference  and the $q$-difference operator in numerical
analysis  and special function.

 \begin{lz}\label{lz25}
Let $f(x,y)=1,g(x,y)=x-y, b_i=x_{i+1},i=0,1,2,\cdots$. Then
\begin{eqnarray}
{\bf\mathbb{D}}^{(n)}_{(g)}[x_1,x_2,\cdots,x_{n+1}]\{F(x)\}=(-1)^nF[x_1,x_2,\cdots,x_{n+1}],
\end{eqnarray}
  where  the classical divided difference $F[x_1,x_2,\cdots,x_{n+1}]$ is recursively defined by
  \begin{eqnarray}
&& F[x_1]=F(x_1); \nonumber\\
   && F[x_1,x_2]=\frac{F(x_1)-F(x_2)}{x_1-x_2}; \nonumber\\
    &&F[x_1,x_2,x_3]=\frac{F[x_1,x_2]-F[x_2,x_3]}{x_1-x_3};\nonumber\\
&&\cdots;\nonumber\\
    &&F[x_1,x_2,\cdots,x_n]=\frac{F[x_1,x_2,\cdots,x_{n-1}]-F[x_2,x_3,\cdots,x_n]}{x_1-x_n}.
\nonumber\end{eqnarray}
  See \cite[p.123]{tttsss} for more details.

As might be expected from this definition by specializing
$b_i=x+h(i-1)$, ${\bf\mathbb{D}}^{(n)}_{(g)}$ can be expressed in a
simpler form
\begin{eqnarray*}
\,\,{\bf\mathbb{D}}^{(n)}_{(g)}[x,x+h,x+2h,\cdots,x+nh]\{F(x)\}=\frac{1}{n!h^n}\nabla^n_hF(x),
 \end{eqnarray*}
where $\nabla_h$ is  the usual  backward  difference operator
defined as
 $\nabla_h\{F(x)\}=F(x)-F(x+h)$. In
the meantime, assume that $F(x)$ is $n$-times differentiable at $x$.
Then it is easily found that
\begin{eqnarray}
\lim_{h\rightarrow
0}\,{\bf\mathbb{D}}^{(n)}_{(g)}[x,x+h,x+2h,\cdots,x+n
h]\{F(x)\}=\frac{1}{n!}\frac{d^n}{dx^n}F(x).
 \end{eqnarray}

Another interesting case is that with the specification
$b_i=xq^i,i=0,1,2,\cdots$, we have
$$\,{\bf\mathbb{D}}^{(1)}_{(g)}[x,xq]\{F(x)\}=\frac{1}{(q-1)}\mathcal{D}_{q,x}F(x),$$
where $\mathcal{D}_{q,x}$ denotes the  $q$-difference operator {\rm
  (cf.\cite{jackson})} appeared previously in (\ref{carbe})
\begin{eqnarray}
 \mathcal{D}_{q,x}F(x)=\frac{F(x)-F(qx)}{x}.\label{dq}
\end{eqnarray}
 \end{lz}Therefore,  from the above observations, we conclude that  the
$n$-th order $(f,g)$-difference operator
$\,{\bf\mathbb{D}}^{(n)}_{(f,g)}$ is particularly useful and
deserves
 further investigation.  As a consequence, we
find that $\,{\bf\mathbb{D}}^{(n)}_{(f,g)}$ satisfies the following
 recursive formula previously well known as a basic property of the classical
divided difference.

\begin{dl}\label{dlnew} Let $f(x,y)\in \mbox{\sl Ker}\mathcal{L}^{(g)}_{3}$, and $\,{\bf\mathbb{D}}^{(n)}_{(f,g)}$ be defined as above.
Then for any  integer $n\geq 0$,
\begin{eqnarray}
\,{\bf\mathbb{D}}^{(n+1)}_{(f,g)}\left[\begin{matrix}b_0,b_{1},\dots ,b_{n+1}\\
x_{1},\dots
,x_{n}\end{matrix}\right]\{F(x)\}&=&\frac{f(x_n,b_{0})}{g(b_{n+1},b_0)}\,{\bf\mathbb{D}}^{(n)}_{(f,g)}\left[\begin{matrix}b_0,b_{1},\dots
,b_{n}\\x_{1},\dots ,x_{n-1}\end{matrix}\right]\{F(x)\}
\nonumber\\
&+&\frac{f(x_n,b_{n+1})}{g(b_0,b_{n+1})}
\,{\bf\mathbb{D}}^{(n)}_{(f,g)}\left[\begin{matrix}b_1,b_{2},\dots
,b_{n+1}\\x_{1},\dots
,x_{n-1}\end{matrix}\right]\{F(x)\}.\label{xinzhi1}
 \end{eqnarray}
 In particular, if $f(x,y)=1\in \mbox{\sl Ker}\mathcal{L}^{(g)}_{3}$, then
\begin{eqnarray}
&&\,{\bf\mathbb{D}}^{(n+1)}_{(g)}[b_0,b_1,b_2,\cdots,b_n,b_{n+1}]\{F(x)\}\\
&=&\frac{\,{\bf\mathbb{D}}^{(n)}_{(g)}[b_0,b_1,b_2,\cdots,b_n]\{F(x)\}-
\,{\bf\mathbb{D}}^{(n)}_{(g)}[b_1,b_2,b_3,\cdots,b_{n+1}]\{F(x)\}}{g(b_{n+1},b_0)}\nonumber.
 \end{eqnarray}
\end{dl}
{\sl Proof.} By Definition \ref{ddd}, a straightforward calculation
leads us to
\begin{eqnarray*}
&&
\,{\bf\mathbb{D}}^{(n+1)}_{(f,g)}\left[\begin{matrix}b_0,b_{1},\dots ,b_{n+1}\\ x_{1},\dots ,x_{n}\end{matrix}\right]\{F(x)\}-\frac{f(x_n,b_{0})}{g(b_{n+1},b_0)}\,{\bf\mathbb{D}}^{(n)}_{(f,g)}\left[\begin{matrix}b_0,b_{1},\dots ,b_{n}\\ x_{1},\dots ,x_{n-1}\end{matrix}\right]\{F(x)\}\\
&=&\sum_{k=0}^{n+1}F(b_k)\frac{\prod_{i=1}^{n} f(x_i,b_k)}
 {\prod_{i=0,i\neq k}^{n+1}g(b_i,b_k)}-\frac{f(x_n,b_{0})}{g(b_{n+1},b_0)}\sum_{k=0}^{n}F(b_k)\frac{\prod_{i=1}^{n-1}
f(x_i,b_k)}
 {\prod_{i=0,i\neq k}^{n}g(b_i,b_k)}\\
 &=&F(b_{n+1})\frac{\prod_{i=1}^{n} f(x_i,b_{n+1})}
 {\prod_{i=0}^{n}g(b_i,b_{n+1})}+\sum_{k=0}^{n}\left\{\frac{f(x_n,b_{k})}{g(b_{n+1},b_k)}-\frac{f(x_n,b_{0})}{g(b_{n+1},b_0)}\right\}
 F(b_k)\frac{\prod_{i=1}^{n-1}
f(x_i,b_k)}
 {\prod_{i=0,i\neq k}^{n}g(b_i,b_k)}.
 \end{eqnarray*}
 Note that the term within the curly braces can be simplified by
 (\ref{idd}), which arises from $f(x,y)\in \mbox{\sl Ker}\mathcal{L}^{(g)}_{3}$. The result is
\begin{eqnarray*}
&&
\sum_{k=0}^{n}\left\{\frac{f(x_n,b_{k})}{g(b_{n+1},b_k)}-\frac{f(x_n,b_{0})}{g(b_{n+1},b_0)}\right\}
 F(b_k)\frac{\prod_{i=1}^{n-1}
f(x_i,b_k)}
 {\prod_{i=0,i\neq k}^{n}g(b_i,b_k)}\\
&=&-\frac{f(x_n,b_{n+1})} {g(b_{n+1},b_0)} \sum_{k=1}^{n}
 F(b_k)\frac{\prod_{i=1}^{n-1}
f(x_i,b_k)}
 {\prod_{i=1,i\neq k}^{n+1}g(b_i,b_k)}.
 \end{eqnarray*}
Simplify the preceding identity by this fact to obtain
\begin{eqnarray*}
&-&\frac{f(x_n,b_{n+1})} {g(b_{n+1},b_0)} \sum_{k=1}^{n}
 F(b_k)\frac{\prod_{i=1}^{n-1}
f(x_i,b_k)}
 {\prod_{i=1,i\neq k}^{n+1}g(b_i,b_k)}+F(b_{n+1})\frac{\prod_{i=1}^{n} f(x_i,b_{n+1})}
 {\prod_{i=0}^{n}g(b_i,b_{n+1})}\\
 &=&\frac{f(x_n,b_{n+1})} {g(b_0,b_{n+1})}
 \sum_{k=1}^{n+1}
 F(b_k)\frac{\prod_{i=1}^{n-1}
f(x_i,b_k)}
 {\prod_{i=1,i\neq k}^{n+1}g(b_i,b_k)}.
 \end{eqnarray*}
By Definition \ref{ddd},  the sum in the r.h.s.\,of this identity is
none other than
$\,{\bf\mathbb{D}}^{(n)}_{(f,g)}\left[\begin{matrix}b_1,b_{2},\dots ,b_{n+1}\\
x_{1},\dots ,x_{n-1}\end{matrix}\right]\{F(x)\}.$ This gives
\begin{eqnarray*}
\,{\bf\mathbb{D}}^{(n+1)}_{(f,g)}\left[\begin{matrix}b_0,b_{1},\dots
,b_{n+1}\\ x_{1},\dots ,x_{n}\end{matrix}\right]\{F(x)\}
&-&\frac{f(x_n,b_{0})}{g(b_{n+1},b_0)}\,{\bf\mathbb{D}}^{(n)}_{(f,g)}\left[\begin{matrix}b_0,b_{1},\dots ,b_{n}\\ x_{1},\dots ,x_{n-1}\end{matrix}\right]\{F(x)\}\\
 &=&\frac{f(x_n,b_{n+1})}{g(b_0,b_{n+1})}
\,{\bf\mathbb{D}}^{(n)}_{(f,g)}\left[\begin{matrix}b_1,b_{2},\dots ,b_{n+1}\\
x_{1},\dots ,x_{n-1}\end{matrix}\right]\{F(x)\}.
\end{eqnarray*}
So the theorem is proved.\qed

A corollary is obtained immediately from the argument of this
theorem.
\begin{tl} Preserve the assumption as above. Then for $n\geq m$,
\begin{eqnarray}
&&{\bf\mathbb{D}}^{(n)}_{(f,g)}\left[\begin{matrix}b_0,b_{1},\dots ,b_{n}\\
x_{1},\dots
,x_{n-1}\end{matrix}\right]\left\{\prod_{i=0}^{m-1}g(b_i,x)/\prod_{i=1}^{m}f(x_i,x)\right\}=
\frac{1}{f(x_m,b_m)}\delta_{n,m};
 \end{eqnarray}
In particular, for\,\,$n\geq 1$,
\begin{eqnarray}
 &&{\bf\mathbb{D}}^{(n)}_{(f,g)}\left[\begin{matrix}b_0,b_{1},\dots ,b_{n}\\
x_{1},\dots ,x_{n-1}\end{matrix}\right]\{1\}=0.
 \end{eqnarray}
\end{tl}
{\sl Proof.}\,\, It can be proved by Theorem \ref{dlnew} and
induction on $n$. Omitted here.\qed

Concerning the $n$th order $(f,g)$-difference of the product of two
functions,  the next theorem states
 that $\,{\bf\mathbb{D}}^{(n)}_{(f,g)}$ also satisfies  the {\sl Leibniz  formula}:

\begin{dl}  Let $f(x,y)\in \mbox{\sl Ker}\mathcal{L}^{(g)}_{3}$, and $\,{\bf\mathbb{D}}^{(n)}_{(f,g)}$ be defined as above. Then
\begin{align}
\,{\bf\mathbb{D}}^{(n)}_{(f,g)}\left[\begin{matrix}b_0,b_{1},\dots
,b_{n}\\ x_{1},\dots
,x_{n-1}\end{matrix}\right]\{F(x)H(x)\}\label{leibniz}\\
=\sum_{k=0}^n
f(x_k,b_k)\,{\bf\mathbb{D}}^{(k)}_{(f,g)}\left[\begin{matrix}b_0,b_{1},\dots
,b_{k}\\ x_{1},\dots
,x_{k-1}\end{matrix}\right]\{H(x)\}\,\,{\bf\mathbb{D}}^{(n-k)}_{(f,g)}\left[\begin{matrix}b_{k},b_{k+1},\dots
,b_{n}\\ x_{k+1},\dots ,x_{n-1}\end{matrix}\right]
\{F(x)\}.\nonumber
\end{align}
\end{dl}
{\sl Proof.} Observe that the r.h.s.\,of  (\ref{leibniz}) can be
reformulated as
\begin{eqnarray}
\sum_{0\leq i\leq j\leq n} \lambda_{i,j}H(b_i)F(b_j),
\end{eqnarray}
where the coefficients $\lambda_{i,j}$ are given by
\begin{eqnarray*}
\lambda_{i,j}&=&\sum_{k=i}^{j}f(x_k,b_k) \frac{\prod_{l=1}^{k-1}
f(x_l,b_i)}
 {\prod_{i_1=0,i_1\neq i}^{k}g(b_{i_1},b_i)}
\frac{\prod_{l=k+1}^{n-1} f(x_l,b_j)}
 {\prod_{i_2=k,i_2\neq j}^{n}g(b_{i_2},b_j)}\\
&=&\frac{\prod_{l=1}^{i-1} f(x_l,b_i)}
 {\prod_{i_1=0}^{i-1}g(b_{i_1},b_i)}
\frac{\prod_{l=j}^{n-1} f(x_l,b_j)}
 {\prod_{i_2=j+1}^{n}g(b_{i_2},b_j)}\\
&\times&\sum_{k=i}^{j}\frac{\prod_{l=i}^{k-1} f(x_l,b_i)}
 {\prod_{i_1=i+1}^{k}g(b_{i_1},b_i)}
\left\{f(x_k,b_k)\frac{\prod_{l=k+1}^{j-1} f(x_l,b_j)}
 {\prod_{i_2=k}^{j-1}g(b_{i_2},b_j)} \right\}.
\end{eqnarray*}
The last sum turns out to be  $\delta_{j,i}$, as the $(j,i)$-th
entry
 of the product of a pair of matrices $(B_{j,k}^{-1})$ and $(B_{k,i})$ given in
Theorem \ref{math4}. It reduces the preceding sum to
\begin{eqnarray}
\sum_{0\leq i\leq j\leq n} \lambda_{i,j}H(b_i)F(b_j)=
\sum_{i=0}^{n}H(b_i)F(b_i)\frac{\prod_{l=1}^{n-1} f(x_l,b_i)}
 {\prod_{i_1=0,i_1\neq i}^{n}g(b_{i_1},b_i)},
\end{eqnarray}
which is just, by Definition \ref{deridef}, the desired result.\qed

Explicit expressions of $\,{\bf\mathbb{D}}^{(n)}_{(f,g)}$ for
$n=1,2$ will now be given as examples to justify this formula.
\begin{lz} Let $\,{\bf\mathbb{D}}^{(n)}_{(f,g)}$  be given by (\ref{deridef}).
Then
\begin{align}
&&\,{\bf\mathbb{D}}^{(1)}_{(g)}[b_0,b_1]\{F(x)H(x)\}=H(b_0)\,{\bf\mathbb{D}}^{(1)}_{(g)}[b_0,b_1]
\{F(x)\}+F(b_1)\,{\bf\mathbb{D}}^{(1)}_{(g)}[b_0,b_1]\{H(x)\};
\label{leib1}
\\
&&\,{\bf\mathbb{D}}^{(2)}_{(f,g)}\left[\begin{matrix}b_0,b_1,b_{2}\\
x_{1}\end{matrix}\right]\{F(x)H(x)\}=f(x_1,b_1)\,{\bf\mathbb{D}}^{(1)}_{(f,g)}[b_0,b_1]\{H(x)\}\,\,{\bf\mathbb{D}}^{(1)}_{(f,g)}[b_1,b_2]\{F(x)\}\nonumber\\
&&+H(b_0)\,{\bf\mathbb{D}}^{(2)}_{(f,g)}\left[\begin{matrix}b_0,b_1,b_{2}\\
x_{1}\end{matrix}\right]\{F(x)\}+F(b_2)\,{\bf\mathbb{D}}^{(2)}_{(f,g)}\left[\begin{matrix}b_0,b_1,b_{2}\\
x_{1}\end{matrix}\right]\{H(x)\}.\label{leib2}
 \end{align}
\end{lz}
We end this section with two connections between
$\,{\bf\mathbb{D}}^{(n)}_{(x-y)}$ and $\mathcal{D}_{q,x}$, as well
as new but short proofs of a few basic propositions of
$\mathcal{D}_{q,x}$. Different from all known arguments is
 that we only utilize the definition of the
$(f,g)$-difference operator.

\begin{xz}\label{xz1}
Let  $\,{\bf\mathbb{D}}^{(n)}_{(g)}$ be given by (\ref{deridef}).
Then the following hold for  $n\geq 1,m\geq 0$,
\begin{eqnarray}
&&\,{\bf\mathbb{D}}^{(n)}_{(x-y)}[x,xq,xq^2,\cdots,xq^n]\{F(x)\}=\frac{(-1)^n}{\poq{q}{n}}\mathcal{D}_{q,x}^n
\left\{F(x)\right\};\label{211}\\
&&\,{\bf\mathbb{D}}^{(n)}_{(x-y)}[xq^m,xq^{m+1},xq^{m+2},\cdots,xq^{n+m}]\{F(x)\}=\frac{(-1)^nq^{-nm}}{\poq{q}{n}}\mathcal{D}_{q,x}^n
\left\{F(xq^m)\right\}.\label{212}
\end{eqnarray}
\end{xz}
{\sl Proof.} We proceed to show (\ref{211}) in question by induction
on $n$. At first, by the definition (see (\ref{dq})), it is easy to
check that
\begin{eqnarray*}
  \mathcal{D}_{q,x}F(x)&=& \frac{F(x)-F(qx)}{x}=(q-1) \frac{F(x)-F(qx)}{x(q-1)}=(q-1)\,{\bf\mathbb{D}}^{(1)}_{(x-y)}[x,xq]\{F(x)\};  \\
 \mathcal{D}^2_{q,x}\{F(x)\}&=&\mathcal{D}_{q,x}\{\mathcal{D}_{q,x}F(x)\}=\frac{qF(x)-(1+q)F(qx)+F(q^2x)}{qx^2}\\
 &=&(q;q)_2\,{\bf\mathbb{D}}^{(2)}_{(x-y)}[x,xq,xq^2]\{F(x)\}.
  \end{eqnarray*}
Hence, the result holds for $n=1,2$. Now, assume that the result in
question is valid for $n=k$. Thus, we have
\begin{eqnarray*}
\,{\bf\mathbb{D}}^{(k)}_{(x-y)}[x,xq,xq^2,\cdots,xq^{k}]\{F(x)\}=
\frac{(-1)^k}{(q;q)_k}\mathcal{D}_{q,x}^k\{F(x)\}.
 \end{eqnarray*}
 Next, consider
$n=k+1$. For this case, from (\ref{xinzhi1}), it follows that
\begin{align*}
\,{\bf\mathbb{D}}^{(k+1)}_{(x-y)}[x,xq,xq^2,\cdots,xq^{k+1}]\{F(x)\}
 =\frac{1}{xq^{k+1}-x}\,{\bf\mathbb{D}}^{(k)}_{(x-y)}[x,xq,xq^2,\cdots,xq^{k}]\{F(x)\}\\
+ \frac{1}{x-xq^{k+1}}
\,{\bf\mathbb{D}}^{(k)}_{(x-y)}[xq,xq^2,\cdots,xq^{k+1}]\{F(x)\}.
 \end{align*}
Observe that
\begin{eqnarray*}
\,{\bf\mathbb{D}}^{(k)}_{(x-y)}[xq,xq^2,\cdots,xq^{k+1}]\{F(x)\}=\,{\bf\mathbb{D}}^{(k)}_{(x-y)}[x,xq,xq^2,\cdots,xq^{k}]\{F(x)\}|_{x\mapsto
qx},
 \end{eqnarray*}
where $x\mapsto qx$ denotes the replacement of $x$ with $xq$. Insert
it into the preceding identity to arrive at
 \begin{eqnarray*}
&& \,{\bf\mathbb{D}}^{(k+1)}_{(x-y)}[x,xq,xq^2,\cdots,xq^{k+1}]\{F(x)\}\\
&=&
\frac{1}{q^{k+1}-1}\frac{\,{\bf\mathbb{D}}^{(k)}_{(x-y)}[x,xq,xq^2,\cdots,xq^{k}]\{F(x)\}-\,{\bf\mathbb{D}}^{(k)}_{(x-y)}
[xq,xq^2,\cdots,xq^{k+1}]\{F(xq)\}}{x}\\
&=&\frac{1}{q^{k+1}-1}\mathcal{D}_{q,x}\left\{\,{\bf\mathbb{D}}^{(k)}_{(x-y)}[x,xq,xq^2,\cdots,xq^{k}]\{F(x)\}\right\}
=\frac{(-1)^{k+1}}{(q;q)_{k+1}}\mathcal{D}_{q,x}\left\{\mathcal{D}^k_{q,x}\{F(x)\}\right\}\\
&=&\frac{(-1)^{k+1}}{\poq{q}{k+1}}\mathcal{D}^{k+1}_{q,x}\{F(x)\}.
 \end{eqnarray*}
Note that the last identity comes from the induction hypothesis.
Thus, the desired result holds also for $n=k+1$. Summing up, this
gives the complete proof of (\ref{211}).

The correctness of (\ref{212}) is proved by noting that
\begin{eqnarray*}
{\bf\mathbb{D}}^{(n)}_{(x-y)}[xq^m,xq^{m+1},xq^{m+2},\cdots,xq^{n+m}]\{F(t)\}\\
=
q^{-nm}\,{\bf\mathbb{D}}^{(n)}_{(x-y)}[x,xq,xq^2,\cdots,xq^n]\{F(tq^m)\}
\end{eqnarray*}
and applying (\ref{211}) to the r.h.s.\,of this identity.
 \qed

Further, the $(f,g)$-difference operator allows to evaluate
$\mathcal{D}^n_{q,x}$ explicitly. As Koornwinder pointed out in his
unpublished \cite{koo}, it was first obtained in 1921 by Ryde
\cite{ryd}.

\begin{xz}\label{yl11}
Let  $\mathcal{D}_{q,x}$ be the usual $q$-difference operator. Then
for $n\geq 1$,
\begin{align}
\mathcal{D}_{q,x}^n\left\{F(x)\right\}=\frac{1}{x^n}\sum_{k=0}^n(-1)^{k}q^{\binom{k+1}{2}-nk}
{n\brack k}_qF(xq^k)\label{272}.
\end{align}
\end{xz}
{\sl Proof.}\,\, Use  Proposition \ref{xz1} to calculate directly
\begin{align*}
\mathcal{D}_{q,x}^n\left\{F(x)\right\}=(-1)^n\poq{q}{n}
\,{\bf\mathbb{D}}^{(n)}_{(g)}[x,xq,xq^2,\cdots,xq^n]\{F(x)\}\\
=(-1)^n\poq{q}{n} \sum_{k=0}^{n}\frac{F(xq^k)}
 {\prod_{i=0,i\neq k}^{n}(xq^i-xq^k)}=\frac{1}{x^n}\sum_{k=0}^n(-1)^{k}q^{\binom{k+1}{2}-nk}
{n\brack k}_qF(xq^k).\end{align*}\qed

By (\ref{272}) and the following finite form of the $q$-binomial
theorem \cite[p.490, Corollary 10.2.2, Eq.(c)]{andrews}
\begin{align}
\sum_{k=0}^{n}(-1)^kq^{\binom{k}{2}}{n\brack
k}_qx^k=\poq{x}{n},\label{binomfinite}
\end{align}
 it is immediately seen that
\begin{xz}\label{xz2} Let $F(x)$ be a polynomial of degree
less than $m$ in $x$. Then for any integer $n\geq m+1$,
\begin{eqnarray}
\mathcal{D}_{q,x}^n\left\{F(x)\right\}=0.
\end{eqnarray}
\end{xz}
\pf Assume that $F(x)=\sum^m_{i=0}a_ix^i$. Then it is easily found
that
\begin{eqnarray*}
\mathcal{D}_{q,x}^n\left\{F(x)\right\}&=&\sum^m_{i=0}a_i\mathcal{D}_{q,x}^n\left\{x^i\right\}=\sum^m_{i=0}a_ix^{i-n}\sum_{k=0}^n(-1)^{k}q^{\binom{k}{2}}
{n\brack k}_qq^{(1+i-n)k}\\
&=&\sum^m_{i=0}a_ix^{i-n}\poq{q^{(1+i-n)}}{n}=0
\end{eqnarray*}
because that for $n\geq m+1\geq i+1\geq 1$,
$\poq{q^{(1+i-n)}}{n}=0$.\qed

 Taken together,  (\ref{leibniz}) and (\ref{272})
leads to
\begin{xz}\label{xz3}{\rm ( See \cite[p.233]{3rr} ).}
\begin{align}
\mathcal{D}_{q,x}^n\left\{F(x)H(x)\right\}=\sum_{k=0}^nq^{(k-n)k}{n\brack
k}_q
\mathcal{D}_{q,x}^k\left\{F(x)\right\}\mathcal{D}_{q,x}^{n-k}\left\{H(q^kx)\right\}.
\end{align}
\end{xz}

\section{The generalized Ismail's argument and  $(f,g)$-expansion formulas}
\setcounter{equation}{0}\setcounter{dl}{0}\setcounter{tl}{0}\setcounter{yl}{0}\setcounter{lz}{0}
Now recall that as a basic application of matrix inversion, if the
two lower-triangular matrices   $(B_{n,k})$ and $(B_{n,k}^{-1})$ are
inverses of each other, then  for any two sequences $\{X_i\}$ and
$\{Y_i\}$,
\begin{equation}\sum _{k=0} ^{n}B_{n,k}X_{k}=Y_{n}\quad \text{if and only if}\quad \sum _{k=0}^{n}
B_{n,k}^{-1}Y_{k}=X_{n}.\label{4.1}\end{equation}
 Simple as it seems, many facts display that (\ref{4.1}) provides a standard and powerful technique for deriving new summation
formulas from known ones. The reader is referred to
\cite{5,6,111,16,kratt, 18,99,milne, schlosser2,1001} for further
details.

In the sequel, once the above standard technique being applied to
Theorem \ref{math4}, it is not hard to set up the following special
result on which  our forthcoming discussions rely.
\begin{yl}\label{yl1}Let $f(x,y)$ and $g(x,y)$, $\{x_i\}$ and $\{b_i\}$ be given as in Theorem
\ref{math4}. Then the system of linear relations for any two
sequences $\{X_i\}$ and $\{Y_i\}$
\begin{eqnarray}X_n=\sum_{k=0}^{n}Y_kf(x_k,b_k)
\frac{\prod_{i=0}^{k-1}g(b_i,b_n)}{\prod_{i=1}^{k}f(x_i,b_n)}\label{2.5}
\end{eqnarray}
 is equivalent to the system
\begin{eqnarray}Y_n=\sum_{k=0}^{n}X_k \frac{\prod_{i=1}^{n-1}f(x_i,b_k)}
 {\prod_{i=0,i\neq k}^{n}g(b_i,b_k)}
 .\label{newc}
\end{eqnarray}
\end{yl}

{\sl Proof.}\,\,Assume that
$$
\sum_{k=0}^{n}X_k \frac{\prod_{i=1}^{n-1}f(x_i,b_k)}
 {\prod_{i=0,i\neq k}^{n}g(b_i,b_k)}
 =Y_n,
$$
or, equivalently
$$
\sum_{k=0}^{n}
 \frac{\prod_{i=k}^{n-1}f(x_i,b_k)}
 {\prod_{i=k+1}^{n}g(b_i,b_k)}\left\{\frac{\prod_{i=1}^{k-1}f(x_i,b_k)}
 {\prod_{i=0}^{k-1}g(b_i,b_k)}X_k\right\} =Y_n.
$$
Then by Theorem \ref{math4}, i.e., the $(f,g)$-inversion, we have
that
$$
\sum_{k=0}^{n}Y_k f(x_k,b_k)\frac{\prod_{i=k+1}^{n-1}f(x_i,b_n)}
{\prod_{i=k}^{n-1}g(b_i,b_n)}=
 \frac{\prod_{i=1}^{n-1}f(x_i,b_n)}
 {\prod_{i=0}^{n-1}g(b_i,b_n)}X_n.
$$
Divide both sides of the relation by $\prod_{i=1}^{n-1}f(x_i,b_n)/
 \prod_{i=0}^{n-1}g(b_i,b_n)$ and simplify the resulted. It gives
$$
X_n=\sum_{k=0}^{n}Y_kf(x_k,b_k)
\frac{\prod_{i=0}^{k-1}g(b_i,b_n)}{\prod_{i=1}^{k}f(x_i,b_n)}.
$$
Vice versa. \qed

If we set $X_n=F(x)|_{x=b_n}$, $F(x)$ is a known analytic function
in certain region, as one of important special case of Lemma
\ref{yl1}, then
$$Y_n={\bf\mathbb{D}}^{(n)}_{(f,g)}\left[\begin{matrix}b_0,b_{1},\dots,b_{n}\\x_{1},\dots,x_{n-1}\end{matrix}\right]\{F(x)\}.
$$
Related by such a pair of sequences, Identity (\ref{2.5}) deserves a
separate definition, which is in analogy to the Fourier series
expansions of analytic functions. For this, we use $\Omega$ to
denote an open subset of the complex plane, and ${\mathcal
H}(\Omega)$ the space of analytic functions over $\Omega$.
\begin{dy}\label{deridef00}With the assumption as above. Let $F(x)\in {\mathcal
H}(\Omega)$. The following series
\begin{eqnarray}\sum_{k=0}^{\infty}G(k)
f(x_k,b_k)
\frac{\prod_{i=0}^{k-1}g(b_i,x)}{\prod_{i=1}^{k}f(x_i,x)}\label{expandef0}
\end{eqnarray}
 is called the
$(f,g)$-series generated by $F(x)$ with respect to two sequences
$\{b_i\}$ and $\{x_i\}$ over $\Omega$.  We denote it by
  \begin{eqnarray}F(x)\sim\sum_{k=0}^{\infty}G(k)
f(x_k,b_k)
\frac{\prod_{i=0}^{k-1}g(b_i,x)}{\prod_{i=1}^{k}f(x_i,x)}\label{expandef1}
\end{eqnarray}
\end{dy}
Very similar to the situation in the theory of Fourier series, we
come up against two questions:
\begin{description}
               \item[The convergence problem] Does the series (\ref{expandef0}) converge at some point $x\in \Omega?$
               \item[The representation problem]If  (\ref{expandef0}) does converge at $x\in \Omega$, is its sum
               $F(x)$? more precisely,
 \begin{eqnarray}F(x)=\sum_{k=0}^{\infty}G(k)
f(x_k,b_k)
\frac{\prod_{i=0}^{k-1}g(b_i,x)}{\prod_{i=1}^{k}f(x_i,x)}.\label{expandef2}
\end{eqnarray}
             \end{description}
In fact, when we deal with these two questions,
 we can not expect a
simple, clear-cut answer without assuming further conditions. In the
present paper, we are mainly  concerned with the second problem with
an effort to set up such expansion formulas. For this, we propose
the following definition
\begin{dy} Let $F(x)\in {\mathcal
H}(\Omega)$ and  $f(x,y)\in \mbox{\sl Ker}\mathcal{L}^{(g)}_{3}$. If
there exist three sequences $\{x_i\},\{b_i\}\subseteq\Omega$, and
$\{\chi(i)\}$, such that for any
$x\in\Omega$,\begin{eqnarray}F(x)=\sum_{k=0}^{\infty}\chi(k)
f(x_k,b_k)
\frac{\prod_{i=0}^{k-1}g(b_i,x)}{\prod_{i=1}^{k}f(x_i,x)}.\label{expandef3}
\end{eqnarray}
then it is called the $(f,g)$-expansion formula of $F(x)$ with
respect to two sequences $\{b_i\}$ and $\{x_i\}$ over $\Omega$.
\end{dy}

The next fact implies that  the $(f,g)$-expansion formula of $F(x)$
is unique if it exists.
\begin{yl} \label{uni} The coefficients $\chi(n)$ in (\ref{expandef3})
are uniquely determined by $\{G(i)\}$, i.e.,
\begin{eqnarray}
\chi(n)=\,{\bf\mathbb{D}}^{(n)}_{(f,g)}\left[\begin{matrix}b_0,b_{1},\dots,b_{n}\\x_{1},\dots,x_{n-1}\end{matrix}\right]\{F(x)\}.
\end{eqnarray}
\end{yl}

It is worth pointing out that even if the $(f,g)$-series generated
by $F(x)$ does not converge to itself, it always agrees with $F(x)$
at infinite points $b_i,i=0,1,\cdots$. Hence, in order to guarantee
that they are equal over $\Omega$, it is sufficient to require that
the series and $F(x)$ be analytic around $x=b$ while $b$ can be
chosen as
 an accumulation point of $\{b_i\}$ in the interior of $\Omega$. Combining
this idea with the original ``Ismail's argument" due to Ismail, we
get all that described in the present paper.

Before proceeding to our main results, we had better display two
examples which can be interpreted in this view. One is Heine's
$q$-analogue of the Gauss summation
    formula (cf.\cite[II.8]{10})
\begin{lz}\label{lz} For four indeterminate $q,a,c,x$ with $|q|<1$ and $|cx/a|<1$,
    \begin{eqnarray}
     && \,_2\phi_1\left[%
\begin{array}{cc}
  a, & 1/x \\
   & c\\
\end{array};q,\frac{cx}{a}
\right]=\frac{(c/a,cx;q)_{\infty}} {(c,cx/a;q)_{\infty}}.\label{68}
    \end{eqnarray}
 \end{lz}
Another is the famous Rogers-Fine identity \cite[p.15, Eq.
(14.1)]{fine}
\begin{lz} \label{finelz}For three indeterminate $q,x,z$ with $|q|<1$ and $|z|<1$,
   \begin{eqnarray}
\sum_{n=0}^{\infty}\frac{\poq{a}{n}}{\poq{x}{n}}z^n=\sum_{k=0}^{\infty}(1-azq^{2k})q^{2\binom{k}{2}}
\frac{\poq{a}{k}(azq/x;q)_k}{(z;q)_{k+1}(x;q)_k}(xz)^k.
\end{eqnarray}
 \end{lz}
Indeed, the former summation formula is an $(1, x-y)$-expansion
formula of the function $F(x)=\frac{(c/a,cx;q)_{\infty}}
{(c,cx/a;q)_{\infty}}$ on the open set $\{x:\, |cx/a|<1\}$ while the
latter transformation formula is an $(1-xy, x-y)$-expansion formula
of the function
$F(x)=\sum_{n=0}^{\infty}\frac{\poq{a}{n}}{\poq{x}{n}}z^n.$

Whereas  many nonterminating  summation and transformation formulas
from the theory of basic
 hypergeometric series, just like Examples \ref{lz} and \ref{finelz}, fit into such a
framework,
 apparently there is no general theorem (aside from those formulas mentioned in Section 1)
 known on the
existence of  such phenomena. Now the definition of
$(f,g)$-expansion formula
 permits us to establish the following version of the
``Ismail's argument".
\begin{dl}[Generalized Ismail's argument]\label{expanth}  Let $F(x)\in {\mathcal
H}(\Omega)$, $f(x,y)\in \mbox{\sl Ker}\mathcal{L}^{(g)}_{3}$. Let
$S_n(x)$ be the sequence of partial sums of the $(f,g)$-series
generated by $F(x)$, say
\begin{eqnarray}S_n(x)=\sum_{k=0}^{n}G(k)
f(x_k,b_k)
\frac{\prod_{i=0}^{k-1}g(b_i,x)}{\prod_{i=1}^{k}f(x_i,x)}\label{expandef1}
\end{eqnarray}
where the coefficients \begin{eqnarray}
G(k)=\,{\bf\mathbb{D}}^{(k)}_{(f,g)}\left[\begin{matrix}b_0,b_{1},\dots,b_{k}\\x_{1},\dots,x_{k-1}\end{matrix}\right]\{F(x)\}.
\label{expandefcoeff1}
\end{eqnarray}
 Assume further that
\begin{description}\item [(i)]  $\lim_{n\mapsto \infty }b_n=b \in
\Omega$;
\item [(ii)] $S_n(x)$  converges uniformly to $S(x)$ in a neighborhood of $b\in
\Omega$;
\item [(iii)] for any integer $i\geq 0, g(b_i,x)/f(x_{i+1},x)$ is analytic
at $x=b$.
   \end{description}
Then there exists  a subset $\Omega_1\subseteq  \Omega$ containing
$b$ such that  for $x\in \Omega_1$,
\begin{eqnarray}
F(x)=S(x)=\sum_{k=0}^{\infty}G(k) f(x_k,b_k)
\frac{\prod_{i=0}^{k-1}g(b_i,x)}{\prod_{i=1}^{k}f(x_i,x)}.\label{2.556}
\end{eqnarray}
\end{dl} {\sl Proof.} At first, we apply  the $(f,g)$-inversion to
the r.h.s.\,of (\ref{expandef1}) to arrive at
\begin{eqnarray}F(b_n)=\sum_{k=0}^{n}G(k)f(x_k,b_k)
\frac{\prod_{i=0}^{k-1}g(b_i,b_n)}{\prod_{i=1}^{k}f(x_i,b_n)},\quad\mbox{thus},\,\,\,
F(b_n)=S_n(b_n).\label{2.55}
\end{eqnarray}
Taken together, the known conditions (ii)-(iii) ensures that $S(x)$
is, by Weierstrass theorem (cf.\cite[p.176,Theorem 1]{ahl}),
analytic at $x=b$. Therefore, taking the limit $n\rightarrow \infty$
on both sides of (\ref{2.55}) and using these basic relations
\begin{eqnarray*}
&&\lim_{n\mapsto\infty} F(b_n)= F(\lim_{n\mapsto\infty}b_n)=F(b);\\
&&S_n(b_n)=\sum_{k=0}^{n}G(k)f(x_k,b_k)
\frac{\prod_{i=0}^{k-1}g(b_i,b_n)}{\prod_{i=1}^{k}f(x_i,b_n)}=\sum_{k=0}^{\infty}G(k)f(x_k,b_k)
\frac{\prod_{i=0}^{k-1}g(b_i,b_n)}{\prod_{i=1}^{k}f(x_i,b_n)};\\
&&\lim_{n\mapsto\infty}S_n(b_n)=\lim_{n\mapsto\infty}\sum_{k=0}^{\infty}G(k)f(x_k,b_k)
\frac{\prod_{i=0}^{k-1}g(b_i,b_n)}{\prod_{i=1}^{k}f(x_i,b_n)}\\
&&=\sum_{k=0}^{\infty}G(k)f(x_k,b_k)
\frac{\prod_{i=0}^{k-1}\lim_{n\mapsto\infty}g(b_i,b_n)}{\prod_{i=1}^{k}\lim_{n\mapsto\infty}f(x_i,b_n)}=S(b),
\end{eqnarray*}
we obtain $F(b)=S(b)$. Based on this, by invoking analytic
continuation argument, there must exist a subset $\Omega_1\subseteq
\Omega$, such that $b\in \Omega_1$ and for $x\in \Omega_1$,
$$ F(x)=\sum_{k=0}^{\infty}G(k) f(x_k,b_k)
\frac{\prod_{i=0}^{k-1}g(b_i,x)}{\prod_{i=1}^{k}f(x_i,x)}=S(x).
$$
So the proof is complete. \qed

\begin{rmm}
 Note that the correctness of (\ref{2.556}) implies that the
set of
    functions
       $\frac{\prod_{i=0}^{k-1}g(b_i,x)}{\prod_{i=1}^{k}f(x_i,x)}$ $(k=0,1,2,\cdots)$
    is a basis of the linear space  ${\mathcal H}(\Omega_1)$ of analytic functions over
$\Omega_1$. For a full characterization of functions $f(x,y)\in
\mbox{\sl Ker}\mathcal{L}^{(g)}_{3}$ such that Condition (iii)
holds, see \cite{0021} for more details.
\end{rmm}
\begin{rmm}\label{counter}
 Note that Conditions {\rm (i)} and {\rm (ii)} are also necessary.  Otherwise, the
conclusion may be false. For instance, let $F(x)=\sin(\pi x)$ and
parameter sequence $b_k=k,k=1,2,\cdots$. Assume that
$$
\sin(\pi x)=\sum_{k=0}^{\infty}G(k)x(x-1)(x-2)\cdots (x-k+1).
$$
By the $(1,x-y)$-inversion, it would follow that
$$G(n)=\,{\bf\mathbb{D}}^{(n)}_{(x-y)}[0,1,\dots,n]\{\sin(\pi\,x)\}\equiv 0$$ for
$n\geq 1$, i.e., $\sin(\pi x)=1$, which is obviously contrary to the
known fact.
\end{rmm}

\begin{rmm}\label{rmm}
As one may expect, the practical difficulty with this theorem
 is the verification of Condition {\rm (ii)}. In this regard,  for the case of $(1,x-y)$-expansion formula, L\'{o}pez-Macro
-Parcet \cite[p.110, Theorem 2.6]{mar} showed that there exists a
disk ${\Bbb O}_r=\{x:|x|<r\}$  such that $S_n(x)$ converges
uniformly with the requirement that $\{b_i\}$ be a $P$-sequence.
Before this, Ismail and Stanton {\rm (cf.\cite[Theorem
3.3]{ismail2})}
 established a $q$-Taylor theorem for any entire function being of  $q$-exponential growth of
 order $2c\ln^2 q$ which is a special case of the
$(1,x-y)$-expansion formula with respect to the sequence $
b_i=(aq^{i}+1/(aq^{i}))/2$. Also, the situation is greatly different
when $F(x)$ is polynomial or  a (basic) hypergeometric series. In
either of these two cases, Condition {\rm (ii)} is considerably easy
to check.
\end{rmm}
Five important cases of Theorem \ref{expanth} are  worthy of note,
which are obtained by specializing  $(f(x,y), g(x,y))=(1,x-y),
(x-y,x-y),(1-xy,x-y)$, and
$((1-axy)(1-b\frac{x}{y}),(x-y)(1-\frac{b}{axy}))$, and
$(y\theta(xy)\theta(\frac{x}{y}),y\theta(xy)\theta(\frac{x}{y}))$,
respectively. We now summarize these results without proof.
\begin{tl} With the same assumption as Theorem \ref{expanth}. Then
the following  hold
 \begin{eqnarray}
\qquad\,F(x)&=&\sum_{k=0}^{\infty}\left\{\sum_{j=0}^{k}\frac{F(b_j)}{\prod_{i=0,i\neq
j}^{k}(b_i-b_j)}\right\}
\prod_{i=0}^{k-1}(b_i-x)\label{1.8}\\
&=&\sum_{k=0}^{\infty}(x_k-b_k)\left\{\sum_{j=0}^{k}F(b_j)\frac{\prod_{i=1}^{k-1}
(x_i-b_j)}
 {\prod_{i=0,i\neq j}^{k}(b_i-b_j)}\right\}
 \frac{\prod_{i=0}^{k-1}(b_i-x)}{\prod_{i=1}^{k}(x_i-x)}\label{1.81}\\
 &=&\sum_{k=0}^{\infty}(1-x_kb_k)\left\{\sum_{j=0}^{k}F(b_j)\frac{\prod_{i=1}^{k-1} (1-x_ib_j)}
 {\prod_{i=0,i\neq j}^{k}(b_i-b_j)}\right\}
 \frac{\prod_{i=0}^{k-1}(b_i-x)}{\prod_{i=1}^{k}(1-xx_i)}\label{1.82}\\
&=&\sum_{k=0}^{\infty}G_1(k)(1-ax_kb_k)(1-bx_k/b_k)
 \frac{\prod_{i=0}^{k-1} (b_i-x)(1-\frac{b}{ab_ix})}{\prod_{i=1}^{k}
 (1-ax_ix)(1-b\frac{x_i}{x})};\label{1.9}\\
&=&\sum_{k=0}^{\infty}G_2(k)b_k\theta(x_kb_k)\theta\left(\frac{x_k}{b_k}\right)
 \frac{\prod_{i=0}^{k-1}\theta(b_ix)\theta(\frac{b_i}{x})}{\prod_{i=1}^{k}
 \theta(x_ix)\theta(\frac{x_i}{x})},\label{1.10}
\end{eqnarray}
where the coefficients
\begin{eqnarray*}
                  G_1(k)&=&\,{\bf\mathbb{D}}^{(k)}_{\left((1-axy)(1-b\frac{x}{y}),(x-y)(1-\frac{b}{axy})\right)}
\left[\begin{matrix}b_0,b_{1},\dots,b_{k}\\x_{1},\dots,x_{k-1}\end{matrix}\right]\{F(x)\},\\
           G_2(k)&=&\,{\bf\mathbb{D}}^{(k)}_{\left(y\theta(xy)\theta(\frac{x}{y}),y\theta(xy)\theta(\frac{x}{y})\right)}
\left[\begin{matrix}b_0,b_{1},\dots,b_{k}\\x_{1},\dots,x_{k-1}\end{matrix}\right]\{F(x)\},\qquad\mbox{and}
                     \end{eqnarray*}
$\theta(x)=(x;q)_\infty(\frac{q}{x};q)_\infty$, $|q|<1$.
\end{tl}
We remark that the theta function $\theta(x)$ satisfying Jacobi's
\emph{triple product identity}
$$
\theta(x)\poq{q}{\infty}=\sum_{k=-\infty}^\infty(-1)^kq^{\binom{k}{2}}x^k,
x\neq 0
$$
 has been frequently used in the study of the elliptic hypergeometric series. See \cite{29,1000} for more details.

 In spite of lacking a general way to verify, as
pointed out in Remark \ref{rmm}, the uniform convergency of
$S_n(x)$,  we have yet the following positive result for the
expansion formula
 (\ref{1.82}).
\begin{dl}[The $(1-xy,x-y)$-expansion theorem]\label{dl32} Let $F(x)\in {\mathcal
H}(\Omega)$, $\{b_i\}, \{x_i\}\subseteq {\Omega}$ such that
$\{b_i\}$ are pairwise distinct and $\{x_i\}$ is bounded,
$\lim_{i\mapsto \infty}b_i=b\neq b_0$ and $\inf\{|1/x_i-b|:i\geq
0\}>0$. Suppose further that $\limsup_{k\mapsto
\infty}|\lambda_k|<\infty$ where
\begin{align}
\lambda_k=\,{\bf\mathbb{D}}^{(k)}_{(1-xy,x-y)}\left[\begin{matrix}b_{1},b_2,\dots,b_{k+1}\\x_{1},
\dots,x_{k-1}\end{matrix}\right]\{F(x)\}/\,{\bf\mathbb{D}}^{(k)}_{(1-xy,x-y)}\left[\begin{matrix}b_0,b_{1},\dots,b_{k}\\x_{1},\dots,x_{k-1}\end{matrix}\right]\{F(x)\}
.\end{align}
 Then  there exists  an open set $\Omega_1$ containing
$b$ such that  for $x\in \Omega_1$,
\begin{align}
F(x)=\sum_{k=0}^{\infty}(1-x_kb_k)\left\{\,{\bf\mathbb{D}}^{(k)}_{(1-xy,x-y)}\left[\begin{matrix}b_0,b_{1},\dots,b_{k}\\x_{1},\dots,x_{k-1}\end{matrix}\right]\{F(x)\}\right\}
 \frac{\prod_{i=0}^{k-1}(b_i-x)}{\prod_{i=1}^{k}(1-xx_i)}.\label{1.822}\end{align}
\end{dl}
{\sl Proof.}\,\, Write $r_0$ for $\inf\{|1/x_i-b|:i\geq 0\}$. The
assumption that $r_0>0$ means that  each function $1/(1-xx_i)$ has
no pole in the disk ${\Bbb O}_{r_0}=\{x:|x-b|<r_0\}$. So according
to Theorem \ref{expanth}, it needs only to check that there exists
an open set $\Omega_1$ containing $b$ , such that $S_n(x)\mapsto
S(x)$ uniformly on $\Omega_1\subseteq {\Bbb O}_{r_0}$. For
notational simplicity, we write
$$
G(k)=\,{\bf\mathbb{D}}^{(k)}_{(1-xy,x-y)}\left[\begin{matrix}b_0,b_{1},\dots,b_{k}\\x_{1},\dots,x_{k-1}
\end{matrix}\right]\{F(x)\}.
$$
 It is clear that by the recursive formula (\ref{xinzhi1}) for the $(f,g)$-difference operator,
the ratio
 \begin{eqnarray*}
\left|\frac{G(k+1)}{G(k)}\right|&=&
\left|\frac{1-x_kb_0}{b_0-b_{k+1}}+\lambda_k\frac{1-x_kb_{k+1}}{b_{k+1}-b_0}\right|.
\end{eqnarray*}
Consider that
 \begin{eqnarray*}
\frac{1-x_kb_0}{b_0-b_{k+1}}+\lambda_k\frac{1-x_kb_{k+1}}{b_{k+1}-b_0}
=x_k+(\lambda_k-1)\frac{1-x_kb_{k+1}}{b_{k+1}-b_0}.\end{eqnarray*}
Thus, by the triangular inequality, we have
 \begin{eqnarray}
\left|\frac{G(k+1)}{G(k)}\right|&\leq&
|x_k|+|\lambda_k-1|\left|\frac{1-x_kb_{k+1}}{b_{k+1}-b_0}\right|.
\end{eqnarray}
 So, given a constant
$M_1=1+\max\{m,(1+|b|m)/|b-b_0|\}$, where $m=\sup\{|x_k|:k\geq 0\}$,
 by solving the inequality under the known condition that
$\lim_{i\mapsto\infty}b_i=b$, we can prove there exists such an
integer $K_1$ that for $k>K$,
$$
\left|\frac{1-x_kb_{k+1}}{b_{k+1}-b_0}\right|<M_1
$$
provided that $M_1\geq \max\{m, (1+m|b|)/|b-b_0|\}$. Next, let
$M=m+m_0M_1$, $m_0=1+\limsup_{k\mapsto \infty}|\lambda_k|$. So for
any $k>K$,
  \begin{eqnarray}
\left|\frac{G(k+1)}{G(k)}\right|\leq
|x_k|+|\lambda_k-1|\left|\frac{1-x_kb_{k+1}}{b_{k+1}-b_0}\right|<m+m_0M_1=M.\label{inequ1}
\end{eqnarray}
Finally,  decide a constant $s<r_0$  by solving the following
inequality for the fixed $M$ and a chosen $a:0<a<1$, $s$ is
independent of $k$ and $K$, such that $|x-b|<s$ small enough, in
order to get
\begin{eqnarray}
&&\left|\frac{1-x_{k+1}b_{k+1}}{1-x_{k}b_{k}}\frac{b_k-x}{1-x_kx}\right|<\frac{a}{M}\left(<\frac{1}{m}\right).\label{inequ2}\end{eqnarray}
 Finally, Inequalities (\ref{inequ1}) and (\ref{inequ1}) together states that for $k>K$,
\begin{eqnarray*}
\left|\frac{\mbox{$(k+1)$-th term  of (\ref{1.822})}}{\mbox{$k$-th
term of
(\ref{1.822})}}\right|=\left|\frac{G(k+1)}{G(k)}\frac{1-x_{k+1}b_{k+1}}{1-x_{k}b_{k}}\frac{b_k-x}{1-x_kx}\right|
<M\times \frac{a}{M}=a<1.\end{eqnarray*}
 By the ratio test, we see that
$S_n(x)$ indeed converges uniformly on the disk ${\Bbb
O}_s=\{x:|x-b|<s\}$, which is such a required open set $\Omega_1$.
Finally, the desired follows from Theorem \ref{expanth}. \qed

Notice that when $F(x)$ is  a basic hypergeometric series, it holds
that $\lim_{k\mapsto \infty}|\lambda_k|=1$. It is also worthy of
note that Fu and Lascoux, in their paper \cite{fu2}, established a
similar expansion formula in the setting of formal power series by
virtue of a divided difference operator acting on
 multivariate function. For a discussion of this divided difference operator and
related matters, the reader may consult \cite{ttt}.
\section{An analytic proof of Gessel and Stanton's $q$-analogue of
the Lagrange inversion formula}Actually, Theorem \ref{dl32} contains
as specifical cases the $q$-analogue of the Lagrange inversion
formula due to Gessel and Stanton (see Lemma \ref{yl112} in Appendix
for a proof of this previously unknown fact), Liu's expansion
formula, and Carlitz's $q$-analogue. We obtain them  only by
specializing $b_i=q^i, x_i=Ap^i$ in (\ref{1.822}). Certainly,  as
initial study on the representation problem of $F(x)$ in terms of
$(f,g)$-series, any possibly rigorous analytic proof of these
formulas is worth investigating. As a result, such a proof is
obtained as follows.
\begin{dl}\label{dl} Let $|p|, |q|<1$ and $F(x)=\sum_{n=0}^\infty a_nx^n$ be arbitrary  power series
with the nonzero radius $R$ of convergence and be also convergent at
$x=R$. Assume further $R$ is not of the form $q^{n} (n\geq 0)$,
$\lim_{n\mapsto\infty}a_{n+1}/a_n=c_0, $ and $0<m |c_0|^n\leq
|a_n|\leq M |c_0|^n$. Then
\begin{eqnarray*}
F(x)&=&\sum_{k=0}^{\infty}G_3(k)(1-Ap^kq^{k})
\frac{\prod_{i=0}^{k-1}(q^i-x)}{\poq{q}{k}(Apx;p)_k}\\
&&\,\,\qquad\qquad\qquad\qquad\quad\mbox{\rm ($\Leftrightarrow$
Gessel and
Stanton's $q$-analogue (\ref{gessel1})/(\ref{gessel2}))}\\
&=&\sum_{k=0}^{\infty}G_4(k)\frac{(1-aq^{2k})(aq/x;q)_kx^k}{(q;q)_k(x;q)_k}\,\,\qquad\mbox{\rm
(Liu's
$q$-expansion formula  (\ref{1111}))}\\
&=&\sum_{k=0}^{\infty}\frac{x^k}{(q,x;q)_k}\left[\mathcal{D}^k_{q,x}\{F(x)(x;q)_{k-1}\}\right]_{x=0}\quad\quad\mbox{\rm
(Carlitz's $q$-analogue
(\ref{carbe}))}\\
\end{eqnarray*}
where the coefficients
\begin{eqnarray}
&&G_3(n)=\sum_{k=0}^n(-1)^{n-k}q^{\binom{k+1}{2}-nk} {n\brack k}_q
(Apq^{k};p)_{n-1}F(q^k);\label{3.17}\\
&&G_4(n)=\frac{1}{a^nq^n}\sum_{k=0}^n(-1)^{k}q^{\binom{k+1}{2}-nk}
\begin{bmatrix}n\\k
\end{bmatrix}_q(aq^{k+1};q)_{n-1}F(aq^{k+1}).\label{271}
\end{eqnarray}
\end{dl}

{\sl Proof.}  Since Liu's $q$-expansion formula (\ref{1111}) follows
from the $q$-analogue of Gessel and Stanton by substituting
$p\mapsto q,A\mapsto a,x\mapsto x/aq$, and Identity (\ref{carbe}) is
the special case $a\mapsto 0$ of (\ref{1111}) together with
Proposition \ref{yl11}, it needs only to show the $q$-analogue of
Gessel and Stanton. To do this, suggested by Theorem \ref{dl32}, it
suffices to check that $\lim_{n\mapsto \infty}|\lambda_n|<\infty$,
because that $b=0,r_0=|1/A|$. Here, it is convenient  to rewrite
$\lambda_n$ in a compact form
\begin{align}
\lambda_n={\Bbb L}_{1:n+1}\{F(x)\}/{\Bbb L}_{0:n}\{F(x)\}\end{align}
with two difference operators employed
\begin{align*}
{\Bbb
L}_{0:n}\{\bullet\}=\,{\bf\mathbb{D}}^{(n)}_{(1-xy,x-y)}\left[\begin{matrix}1,q,q^2,\dots,q^{n}\\Ap,Ap^2,\dots,Ap^{n-1}
\end{matrix}\right]\{\bullet\};\\
{\Bbb
L}_{1:n+1}\{\bullet\}=\,{\bf\mathbb{D}}^{(n)}_{(1-xy,x-y)}\left[\begin{matrix}
q,q^2,q^3,\dots,q^{n+1}\\Ap,Ap^2,\dots,Ap^{n-1}\end{matrix}\right]\{\bullet\}.\end{align*}
Our consideration breaks up into the following  two cases.\\

 {\bf Case I.} At first, let $F(x)=x^r$. By using of the fact
\begin{align}
 \sum_{k=0}^{n}
 \frac{b_k^{r}}{\prod_{i=0,i\neq k}^{n}(b_i-b_k)}=(-1)^n h_{r-n}(b_0,b_1,b_2,\cdots,b_n),
\end{align}
it is easily found that
\begin{align}
\,{\bf\mathbb{D}}^{(n)}_{(1-xy,x-y)}\left[\begin{matrix}b_0,b_1,b_{2},\cdots,b_{n}\\
x_{1},x_2,\cdots,x_{n-1}\end{matrix}\right]\{x^r\}
=\sum_{k=0}^{n-1}(-1)^{n-k}e_{k}(x_{1},x_2,\cdots,x_{n-1})h_{r+k-n}(b_0,b_1,b_2,\cdots,b_n),\nonumber
\end{align}
where $h_i$ (resp. $e_i$), as the $i$th complete symmetric function
of $b_0,b_1,\cdots,b_n$ (resp. $x_1,x_2,\cdots,x_{n-1}$), is given
by
\begin{eqnarray*}
&&h_r(b_0,b_1,b_2,\cdots,b_n)=\sum_{0\leq i_1\leq i_2\leq \cdots\leq
i_r\leq n}b_{i_1}b_{i_2}\cdots b_{i_r};\label{complete}
\\
&&e_{k}(x_1,x_2,\cdots,x_{n-1})=\sum_{1\leq i_1< i_2<\cdots< i_k\leq
n-1}x_{i_1}x_{i_2}\cdots x_{i_k}.\label{complete1}
\end{eqnarray*}
With the choices that $b_i=q^i,x_i=Ap^i$, we establish by induction
on $r$ (resp. $k$) that
\begin{eqnarray*} &&h_r(1,q,q^2,\cdots,q^n)=\sum_{0\leq i_1\leq
i_2\leq \cdots\leq i_r\leq n}q^{i_1+i_2+\cdots+i_r}={n+r\brack n}_q;
\\
&&e_{k}(Ap,Ap^2,\cdots,Ap^{n-1})=\sum_{1\leq i_1< i_2<\cdots<
i_k\leq
n-1}A^kp^{i_1+i_2+\cdots+i_k}=A^kp^{\binom{k+1}{2}}{n-1\brack k}_p.
\end{eqnarray*}
These results combined together leads us to
\begin{eqnarray*}
{\Bbb L}_{0:n}\{x^r\}&=&
\sum_{k=0}^{n-1}(-1)^{n-k}A^{k}p^{\binom{k+1}{2}} {n-1\brack k}_p
{r+k\brack n}_q\\
&=&\sum_{k=1}^{r}(-1)^{k}A^{n-k}p^{\binom{n-k+1}{2}} {n-1\brack
n-k}_p {r+n-k\brack r-k}_q;
\end{eqnarray*}
\begin{eqnarray*}
{\Bbb L}_{1:n+1}\{x^r\}&=&
\sum_{k=0}^{n-1}(-1)^{n-k}A^{k}p^{\binom{k+1}{2}} {n-1\brack k}_p
q^{r+k-n} {r+k\brack
n}_q\\&=&\sum_{k=1}^{r}(-1)^{k}A^{n-k}p^{\binom{n-k+1}{2}}q^{r-k}
{n-1\brack n-k}_p {r+n-k\brack r-k}_q.
\end{eqnarray*}
Finally, we calculate directly
\begin{align}
\lim_{n\mapsto \infty}\lambda_n=\lim_{n\mapsto \infty}\frac{{\Bbb
L}_{1:n+1}\{x^r\}}{{\Bbb L}_{0:n}\{x^r\}}=\lim_{n\mapsto
\infty}\frac{\sum_{k=1}^{r}(-1)^{n-k}A^{n-k}p^{\binom{n-k+1}{2}}q^{r-k}
{n-1\brack k-1}_p {r+n-k\brack
r-k}_q}{\sum_{k=1}^{r}(-1)^{n-k}A^{n-k}p^{\binom{n-k+1}{2}}
{n-1\brack n-k}_p {r+n-k\brack
r-k}_q}\nonumber\\
=\frac{\sum_{k=1}^{r}\lim_{n\mapsto
\infty}(-1)^{k}A^{-k}p^{n(r-k)+\binom{k}{2}}q^{r-k} {n-1\brack
k-1}_p {r+n-k\brack r-k}_q}{\sum_{k=1}^{r}\lim_{n\mapsto
\infty}(-1)^{k}A^{-k}p^{n(r-k)+\binom{k}{2}} {n-1\brack k-1}_p
{r+n-k\brack r-k}_q}=1
\end{align}
by noting that
$$
\lim_{n\mapsto \infty}(-1)^{k}A^{-k}p^{n(r-k)+\binom{k}{2}}
{n-1\brack k-1}_p {r+n-k\brack r-k}_q= \left\{\begin{array}{ll}
             0, & k\neq r;\\
             (-1)^{r}A^{-r}p^{\binom{r}{2}}
\frac{1}{\poqp{p}{r-1}},&  k=r.
           \end{array}
         \right.
$$
  {\bf Case II.} Next, suppose that $F(x)=\sum_{r=0}^{\infty}a_rx^r$ is a power series with the nonzero radius $R$
  of convergency. By the definition, we have
\begin{eqnarray}
{\Bbb L}_{1:n+1}\{F(x)\}&=&
\sum_{k=1}^{n}(-1)^{k}A^{n-k}p^{\binom{n-k+1}{2}} {n-1\brack
n-k}_p\sum_{r=k}^{\infty}a_rq^{r-k} {r+n-k\brack
r-k}_q\nonumber\\
&=&\sum_{k=1}^{n}(-1)^{k}A^{n-k}p^{\binom{n-k+1}{2}} {n-1\brack
n-k}_p\sum_{r=0}^{\infty}a_{r+k}q^{r} {r+n\brack r}_q.\end{eqnarray}
For the sake of simplicity,  define that
\begin{eqnarray}
K_{n,k}(x)&=&\sum_{r=0}^{\infty}a_{r+k}{r+n\brack
r}_qx^{r}.\label{knk}
\end{eqnarray}
By this definition, it is easily found that the generating function
of  $K_{n,k}(x)$ in $t: |t|<1$
\begin{eqnarray}
&&\sum_{k=0}^{\infty}K_{n,k}(x)t^k=\frac{F(t)}{\poq{x/t}{n+1}}\qquad\mbox{and}\label{4.8}\\
&&\mathcal{D}^n_{q,x}\{F(x)\}=\poq{q}{n}K_{n,n}(x).
\end{eqnarray}
Further, from the known condition $\lim_{n\mapsto\infty}
a_{n+1}/a_n=c_0$, it follows
\begin{align}
\lim_{n\mapsto\infty}\frac{K_{n,n}(x)}{a_n}=
\sum_{r=0}^{\infty}\lim_{n\mapsto\infty}\frac{a_{r+n}}{a_n}
{r+n\brack r}_q x^{r}=\sum_{r=0}^{\infty}
\frac{(xc_0)^r}{\poq{q}{r}}=\frac{1}{\poq{xc_0}{\infty}}.\label{result}
\end{align}
Note that, by Hadamard's formula, $|c_0|=1/R$. Two recursive
relations of interest implied by (\ref{4.8}) are
\begin{eqnarray*}
&&K_{n,k}(x)-K_{n,k}(qx)=xK_{n,k+1}(x)-xq^{n+1}K_{n,k+1}(qx);\\
&&K_{n,k}(x)=\frac{1}{1-q^n}K_{n-1,k}(x)-\frac{q^{n}}{1-q^n}K_{n-1,k}(qx).
\end{eqnarray*}
By  iterating the last relation $m$ times, we obtain at once
\begin{eqnarray}
K_{n,k}(x) &=&\frac{\poq{q}{n-m}}{\poq{q}{n}}
\sum_{i=0}^m(-1)^iq^{(n-m+1)i+\binom{i}{2}}{m\brack
i}_qK_{n-m,k}(q^ix).\label{3.22}
\end{eqnarray}
Next, replacing $k$ by $n-k$ and  letting $m=k$ simultaneously in
(\ref{3.22}) gives
\begin{eqnarray*}
K_{n,n-k}(x)&=&\frac{\poq{q}{n-k}}{\poq{q}{n}}
\sum_{i=0}^k(-1)^iq^{(n-k+1)i+\binom{i}{2}}{k\brack
i}_qK_{n-k,n-k}(q^ix).
\end{eqnarray*}
Reformulate by this fact the sum
\begin{eqnarray*}
&&\sum_{k=1}^{n}(-1)^{k}A^{n-k}p^{\binom{n-k+1}{2}} {n-1\brack
n-k}_pK_{n,k}(x)=\sum_{k=0}^{n-1}(-1)^{n-k}A^{k}p^{\binom{k+1}{2}}
{n-1\brack k}_pK_{n,n-k}(x)\\
&=&\sum_{k=0}^{n-1}(-1)^{n-k}A^{k}p^{\binom{k+1}{2}} {n-1\brack
k}_p\frac{\poq{q}{n-k}}{\poq{q}{n}}
\sum_{i=0}^k(-1)^iq^{(n-k+1)i+\binom{i}{2}}{k\brack
i}_qK_{n-k,n-k}(q^ix).
\end{eqnarray*}
Now, we turn to the limitation of ${\Bbb L}_{1:n+1}\{F(x)\}/{\Bbb
L}_{0:n}\{F(x)\}$ as $n\mapsto\infty$. For this we may split ${\Bbb
L}_{1:n+1}\{F(x)\}$ and ${\Bbb L}_{0:n}\{F(x)\}$ into
\begin{eqnarray*} &&{\Bbb
L}_{1:n+1}\{F(x)\}=(-1)^nK_{n,n}(q)+S(n,q)\qquad\qquad\mbox{and}\\
&&{\Bbb L}_{0:n}\{F(x)\}=(-1)^nK_{n,n}(1)+S(n,1)
\end{eqnarray*}
by defining
\begin{eqnarray*}
S(n,x)&=&\sum_{k=1}^{n-1}(-1)^{n-k}A^{k}p^{\binom{k+1}{2}}
{n-1\brack k}_p\frac{\poq{q}{n-k}}{\poq{q}{n}}
\sum_{i=0}^k(-1)^iq^{(n-k+1)i+\binom{i}{2}}{k\brack
i}_qK_{n-k,n-k}(q^ix).
\end{eqnarray*}
By using of  Tannery's theorem\,(cf.\cite{boa} or \cite[Appendix:
Tannery's limiting theorem]{chu}), we obtain that
 $\lim_{n\mapsto\infty}S(n,x)/a_n=0$ \footnote{For its proof, see Lemma \ref{bu2} in
Appendix.} while $x=1,q$. This leads us to
\begin{eqnarray*}
\lim_{n\mapsto\infty}|\lambda_n|&=&\lim_{n\mapsto\infty}\left|\frac{{\Bbb
L}_{1:n+1}\{F(x)\}}{{\Bbb
L}_{0:n}\{F(x)\}}\right|=\lim_{n\mapsto\infty}\left|\frac{(-1)^nK_{n,n}(q)/a_n+S(n,q)/a_n}{(-1)^nK_{n,n}(1)/a_n
+S(n,1)/a_n}\right|\\
&=&\lim_{n\mapsto\infty}\left|\frac{(-1)^nK_{n,n}(q)/a_n}{(-1)^nK_{n,n}(1)/a_n}\right|=|\frac{\poq{c_0}{\infty}}{\poq{qc_0}{\infty}}|=|1-c_0|\leq
1+1/R.\label{d0}
\end{eqnarray*}
The latter arises from (\ref{result}).  This gives the complete
proof of theorem.\qed

 This corollary states clearly that, just as mentioned earlier, the $(1-xy,x-y)$-expansion formula (\ref{1.822})
indeed unifies Gessel and Stanton's $q$-analogue,  Liu's expansion
formula, and Carlitz's $q$-analogue  of the Lagrange inversion
formula. The approach we have presented as above is different
 from  that of Liu \cite{liu}, who  established (\ref{1111}) by means of (\ref{carbe}) and
   expressed  the coefficients $G_2(k)$ in terms of $\mathcal{D}^i_{q,x},i=0,1,2,...$, which is also different from  (\ref{271}).

 As an illustrative example, the calculation of ${\Bbb
L}_{0:n}\{\frac{1}{1-c x}\}$ ( $c\neq 0$) leads us to an interesting
expansion formula.
\begin{lz}
Let $F(x)=\frac{1}{1-c
 x},\, |cx|<1,c\neq 0$. In this case, we calculate by the definition
 \begin{eqnarray*}
{\Bbb L}_{0:n}\{\frac{1}{1-c x}\}&=&\sum_{r=0}^{\infty}c^r{\Bbb
L}_{0:n}\{x^r\}=\sum_{r=0}^{\infty}c^r\sum_{k=1}^{r}(-1)^{k}A^{n-k}p^{\binom{n-k+1}{2}}
{n-1\brack n-k}_p {r+n-k\brack r-k}_q\\&=&
\sum_{k=1}^{n}(-1)^{k}A^{n-k}p^{\binom{n-k+1}{2}} {n-1\brack
n-k}_p\sum_{r=k}^{\infty}c^r {r+n-k\brack r-k}_q\\
&=& \frac{(-1)^nc^n}{\poq{c}{n+1}}
\sum_{k=0}^{n-1}(-1)^{k}(Ap/c)^{k}p^{\binom{k}{2}} {n-1\brack k}_p
=(-1)^n\frac{c^n\poqp{Ap/c}{n-1}}{\poq{c}{n+1}}. \end{eqnarray*}
Similarly, we have
\begin{eqnarray*}
{\Bbb L}_{1:n+1}\{\frac{1}{1-c x}\} &=&
\frac{(-1)^nc^n}{\poq{cq}{n+1}}
\sum_{k=0}^{n-1}(-1)^{k}(Ap/c)^{k}p^{\binom{k}{2}} {n-1\brack k}_p
=(-1)^n\frac{c^n\poqp{Ap/c}{n-1}}{\poq{cq}{n+1}}.
\end{eqnarray*}
This not only gives that $\lim_{n\mapsto \infty}\lambda_n=1-c$ but
also
 leads to
 the expansion formula
\begin{eqnarray}
\frac{1}{1-c
x}&=&\frac{1}{(1-c)(1-A/c)}\sum_{k=0}^{\infty}\frac{1-Ap^kq^{k}}{\poq{q}{k}}\frac{\poqp{A/c}{k}}{\poq{cq}{k}}
\frac{\poq{1/x}{k}}{(Apx;p)_k}\,(cx)^k.
\end{eqnarray}
The special case  $A=0$ of this result gives a refined formula
\begin{eqnarray} \frac{1-c}{1-c
x}&=&\sum_{k=0}^{\infty}\frac{\poq{1/x}{k}}{\poq{q,cq}{k}} (cx)^k.
\end{eqnarray}
Its equivalent form given by (\ref{3.17}) is
\begin{eqnarray}
\frac{c^n}{\poq{c}{n+1}}&=&\frac{1}{\poq{q}{n}}\sum_{k=0}^{n}(-1)^{k}q^{\binom{k+1}{2}-nk}{n\brack
k}_q\frac{1}{1-c q^k}.
\end{eqnarray}
Restate it in terms of $\mathcal{D}_{q,x}$. The result is
\begin{eqnarray}
\mathcal{D}^n_{q,x}\left\{\frac{1}{1-c
x}\right\}|_{x=1}=\frac{c^n\poq{q}{n}}{\poq{c}{n+1}}.
\end{eqnarray}
\end{lz}
\section{ Applications to basic hypergeometric series}\setcounter{equation}{0}\setcounter{dl}{0}\setcounter{tl}{0}\setcounter{yl}{0}\setcounter{lz}{0}

Our applications of all results in the forgoing sections to basic
hypergeometric series is realized two-fold way. First, start with a
given terminating summation and a properly chosen parameter sequence
$b_n\mapsto b$, we want to extend it to a nonterminating summation
as an $(f,g)$-expansion formula of certain analytic function,
proceeding as the ``Ismail's argument", by analytic continuity on
the given region.  Second, just the reverse of the above procedure,
for any known $(f,g)$-expansion of $F(x)$, there always exists a
terminating summation formula, i.e., the expression of
 the $n$-th order $(f,g)$-difference $G(n)$ of $F(x)$ in terms of $F(x)|_{x=b_i}$. For the limitation
of space, we only consider a few remarkable summation and
transformation formulas.
\subsection{From terminating summations to nonterminating summations}Since it is the origin of the ``Ismail's argument", we
would like to adopt Ismail's proof of {\rm Ramanujan's
${}_{1}\varphi_{1}$ summation formula as our first example. See
\cite{ismail0} or \cite[pp.504-505]{andrews} for more details.

\begin{dl} For $|q|<1$ and
$|ba^{-1}|<|x|<1$,
    \begin{eqnarray}
     &&
\sum_{n=-\infty}^{\infty}\frac{\poq{a}{n}}{\poq{b}{n}}x^n=\frac{(q,b/a,ax,q/(ax);q)_{\infty}}
{(b,q/a,x,b/(ax);q)_{\infty}}.
    \end{eqnarray}
      \end{dl}
{\sl Proof.}  At first,  define two functions of a variable $y$
\begin{eqnarray*}
 S(y)&=&
 \sum_{n=-\infty}^{\infty}\frac{\poq{a}{n}}{\poq{y}{n}}x^n\qquad\mbox{and}\quad F(y)=\frac{(q,y/a,ax,q/(ax);q)_{\infty}}
{(y,q/a,x,y/(ax);q)_{\infty}}.
    \end{eqnarray*}
 Evidently, $F(y)$ and $S(y)$ are analytic around $y=0=\lim_{N\mapsto\infty}
b_N$, $b_N=q^{N+1}, N\geq 0$. And then to show $F(b_N)=S(b_N)$ for
each $N$, i.e.,
\begin{eqnarray}
     &&
\sum_{n=-\infty}^{\infty}\frac{\poq{a}{n}}{\poq{q^{N+1}}{n}}x^n=\frac{(q,q^{N+1}/a,ax,q/(ax);q)_{\infty}}
{(q^{N+1},q/a,x,q^{N+1}/(ax);q)_{\infty}},
    \end{eqnarray}
which follows from rewriting the $q$-binomial theorem as follows
\begin{eqnarray}
     && \frac{\poq{ax}{\infty}}{\poq{x}{\infty}}=\sum_{n=0}^{\infty}\frac{\poq{a}{n}}{\poq{q}{n}}x^n
=\sum_{n=-N}^{\infty}\frac{\poq{a}{n+N}}{\poq{q}{n+N}}x^{n+N}
\nonumber\\&=&\sum_{n=-N}^{\infty}\frac{\poq{a}{n+N}}{\poq{q}{n+N}}x^{n+N}=\frac{(a;q)_{N}}{(q;q)_{N}}x^N
\sum_{n=-N}^{\infty}\frac{\poq{aq^N}{n}}{\poq{q^{N+1}}{n}}x^{n}
    \end{eqnarray}
and then replacing $a$ with $aq^{-N}$ in the above identity. Note
that by the convention or the inverse identity (cf.\cite[I.2]{10})
$$
\frac{1}{\poq{q^{N+1}}{n}}=0\qquad\mbox{for}\quad n<-N-1.
$$
By analytic continuation, we see that the result is true. \hfill\qed

Another example of the $(1-xy,x-y)$-expansion formula, also a good
 case to illustrate the ``Ismail's argument", is the following generalized Lebesgue identity due to
Carlitz \cite{carliz1}.
\begin{dl}
\label{finalz}For three indeterminate $q,x,a$ with $|q|<1$ and
$|x|<1$,
   \begin{eqnarray}
\sum_{k=0}^{\infty}(-1)^kq^{\binom{k}{2}}
\frac{\poq{x}{k}a^k}{(q,bx;q)_k}=\frac{\poq{a,x}{\infty}}{\poq{bx}{\infty}}
\sum_{k=0}^{\infty}\frac{\poq{b}{k}}{\poq{q,a}{k}}x^k.
\end{eqnarray}
\end{dl}
\pf At first, apply  the substitution $b\mapsto b/y,x\mapsto
xy,a\mapsto y$ to reformulate equivalently the identity in question
as
\begin{eqnarray}
F(y)=\sum_{k=0}^{\infty}\frac{(-1)^kx^k}{\poq{q}{k}}
\frac{\prod_{i=0}^{k-1}g(b_i,y)}{\prod_{i=1}^{k}f(x_i,y)},\end{eqnarray}
where $f(x,y)=1-xy,g(x,y)=x-y, b_i=bq^i,x_i=q^{i-1}$, and $F(y)$ is
defined by
\begin{eqnarray}
F(y)=\sum_{i=0}^{\infty}\frac{(-1)^iq^{\binom{i}{2}}\poq{bxq^i}{\infty}}{\poq{q}{i}}\frac{y^i}{\poq{y,xyq^i}{\infty}}.
\end{eqnarray}
By the $(1-xy,x-y)$-expansion formula, it suffices to show in a
direct way that
\begin{eqnarray}
&&\frac{(bx)^n}{1-bq^{2n-1}}=\sum_{k=0}^n(-1)^{k}q^{\binom{k+1}{2}-nk}
{n\brack k}_q (bq^{k};q)_{n-1}F(bq^k).\label{33.17}
\end{eqnarray}
For this, we calculate by using Heine's transformation formula to
obtain
\begin{eqnarray*}
F(bq^k)&=&\frac{\poq{bx}{k}}{\poq{bq^k}{\infty}}
\sum_{i=0}^{\infty}\frac{(-1)^iq^{\binom{i}{2}}(bq^k)^i}{\poq{q}{i}}\frac{\poq{bxq^k}{i}}{\poq{bx}{i}}\\
&=&\sum_{i=0}^{\infty}\frac{\poq{q^{-k}}{i}(bxq^k)^i}{\poq{q,bq^k}{i}}=\sum_{i=0}^{k}(-1)^iq^{\binom{i}{2}}{k\brack
i}_q \frac{(bx)^i}{\poq{bq^k}{i}}.
\end{eqnarray*}
Inserting this into the r.h.s. of (\ref{33.17}) to simplify the
resulting identity, we have
\begin{eqnarray*}
\mbox{RHS of (\ref{33.17})}&=&\sum_{n\geq k\geq i\geq
0}(-1)^{k+i}q^{\binom{k+1}{2}-nk+\binom{i}{2}}{n\brack k}_q{k\brack
i}_q\frac{\poq{bq^k}{n-1}(bx)^i}{\poq{bq^k}{i}}\\
&=&\sum_{i=0}^n(bx)^iq^{-Ni}{n\brack
i}_q\sum_{K=0}^N(-1)^{K}q^{\binom{K+1}{2}-NK}{N\brack
K}_q\frac{\poq{b}{N+K+2i-1}}{\poq{b}{K+2i}},
\end{eqnarray*}
where $N=n-i, K=k-i$. Evidently, if $N\geq 1$, then
$$
\frac{\poq{b}{N+K+2i-1}}{\poq{b}{K+2i}}=\poq{bq^{K+2i}}{N-1}.
$$
Since for fixed $i$, $\poq{bq^{2i}x}{N-1}$ is a polynomial of degree
no more than $N-1$ in $x$, thus by Proposition \ref{xz2}, we find
that the inner sum turns out to be
\begin{eqnarray*}
\sum_{K=0}^N(-1)^{K}q^{\binom{K+1}{2}-NK}{N\brack
K}_q\frac{\poq{b}{N+K+2i-1}}{\poq{b}{K+2i}}=0.
\end{eqnarray*}
Thus, the only nonzero term corresponding to  $N=0$ of the r.h.s. of
(\ref{33.17}) is $(bx)^n/(1-bq^{2n-1})$.
 It gives the complete proof of (\ref{33.17}).
\qed

 Next, by the $(f,g)$-inversion and the
$q$-Pfaff-Saalsch\"{u}tz formula (which makes our proofs given here
different from all known ones), together with the ``Ismail's
argument", we obtain new proofs of Rogers' nonterminating
very-well-poised (in short, VWP-balanced) ${}_{6}\phi _{5}$
summation formula (cf.\cite[Eq.(II.20)]{10}) and Watson's
   $q$-analogue of Whipple's transformation formula from
  ${}_{8}\phi _{7}$ to ${}_{4}\phi _{3}$ (cf.\cite[III.17]{10}).
\begin{dl}[Rogers' nonterminating VWP-balanced ${}_{6}\phi _{5}$ summation
formula] For all parameters $a,b,c,d$ where the series converges,
    \begin{align}
     \,_6\phi_5\left[%
\begin{array}{cccccccc}
  a, & q\sqrt{a}, & -q\sqrt{a}, & b, & c, & d \\
   & \sqrt{a}, & -\sqrt{a}, & aq/b, & aq/c, & aq/d  \\
\end{array};q,\frac{aq}{bcd}
\right]=\frac{(aq,aq/(cd),aq/(bd),aq/(bc);q)_{\infty}}
{(aq/b,aq/c,aq/d,aq/(bcd);q)_{\infty}}.\label{688}
    \end{align}
\end{dl}
 {\sl Proof.} Note that the $q$-Pfaff-Saalsch\"{u}tz formula\begin{eqnarray}
\,_3\phi_2\left[%
\begin{array}{cccccccc}
  q^{-n}, & aq^n, & aq/bc \\
   & aq/b, & aq/c \\
\end{array};q,q
\right]=\left(\frac{aq}{bc}\right)^n\frac{\poq{b,c}{n}}{\poq{aq/b,aq/c}{n}}
\end{eqnarray}
 can be reformulated as
\begin{align}&&\sum_{k=0}^n\frac{\poq{q^{-n}}{k}}{\poq{q}{k}}\frac{\poq{aq^{n}}{k}}{\poq{aq}{k}}q^k
 \frac{\poq{aq,aq/(bc)}{k}}
{\poq{aq/b,aq/c}{k}}=\frac{\poq {b,c}{n}} {\poq
{aq/b,aq/c}{n}}\left(\frac{aq}{bc}\right)^n.\label{saaa}\end{align}
 Here,
we will invoke the following  $(x-y,x-y)$-inversion originally due
to Carlitz  (cf.\cite{6} or \cite[Corollary 1]{0022})
\begin{eqnarray*}
\left(\frac{\poq{q^{-n}}{k}}{\poq{q}{k}}\frac{\poq{aq^{n}}{k}}{\poq{aq}{k}}q^k\right)^{-1}=\left(\frac{\poq{q^{-n}}{k}}{\poq{aq^{1+n}}{k}}\frac{\poq{a}{k}}{\poq{q}{k}}\frac{1-aq^{2k}}{1-a}q^{kn}\right).
 \end{eqnarray*}
Apply this inversion to (\ref{saaa}) to get
\begin{align}
 \,_6\phi_5\left[%
\begin{array}{cccccccc}
  a, & q\sqrt{a}, & -q\sqrt{a}, & b, & c, & q^{-n} \\
   & \sqrt{a}, & -\sqrt{a}, & aq/b, & aq/c, & aq^{1+n}  \\
\end{array};q,\frac{aq^{1+n}}{bc}
\right]=\frac{\poq{aq,aq/(bc)}{n}} {\poq{aq/b,aq/c}{n}}.\label{nnn}
    \end{align} Define further
\begin{eqnarray*}
 S(x)&=&
 \,_6\phi_5\left[%
\begin{array}{cccccccc}
  a, & q\sqrt{a}, & -q\sqrt{a}, & b, & c, & 1/x \\
   & \sqrt{a}, & -\sqrt{a}, & aq/b, & aq/c, & aqx  \\
\end{array};q,\frac{aqx}{bc}
\right];\\
F(x)&=&\frac{(aq,aq/(bc);q)_{\infty}}
{(aq/b,aq/c;q)_{\infty}}\frac{(aqx/c,aqx/b;q)_{\infty}}
{(aqx,aqx/(bc);q)_{\infty}}.
    \end{eqnarray*}
 Evidently, $F(x)$ and $S(x)$ are analytic around $x=0$. Choose that $b_n=q^n, x_n=aq^n$. Note that $f(x_n,b_n)=1-aq^{2n}, \lim_{n\mapsto \infty}b_n=0$,
and $ F(b_n)=S(b_n)$. Finally, by Theorem \ref{dl}, the claimed
 follows.
  \hfill\qed

 \begin{dl}[Watson's
   $q$-analogue of Whipple's transformation formula from
  ${}_{8}\phi _{7}$ to ${}_{4}\phi _{3}$]
       \begin{align}
&&    \,_8\phi_7\left[%
\begin{array}{cccccccc}
  a, & q\sqrt{a}, & -q\sqrt{a}, & b, & c, & d,&e,&f \\
   & \sqrt{a}, & -\sqrt{a}, & aq/b, & aq/c, & aq/d, & aq/e, & aq/f\\
\end{array};q,\frac{a^2q^2}{bcdef}
\right]\label{waston}\\
&&=\frac{(aq,aq/(ef),aq/(df),aq/(de);q)_{\infty}}
{(aq/d,aq/e,aq/f,aq/(def);q)_{\infty}}\,_4\phi_3\left[%
\begin{array}{cccccccc}
  aq/(bc), & d, & e, & f \\
   & def/a ,& aq/b, & aq/c \\
\end{array};q,q
\right]\nonumber
    \end{align}
provided the $\,_4\phi_3$ series terminates.
\end{dl}
{\sl Proof.}
    By the $q$-Pfaff-Saalsch\"{u}tz formula, we obtain
\begin{eqnarray}
&&\,_3\phi_2\left[%
\begin{array}{cccc}
  q^{-n}, & aq^{n}, & aq/bc\\
   & aq/b, & aq/c \\
\end{array};q,q
\right]=\left(\frac{aq}{bc}\right)^{n}\frac{\poq{b,c}{n}}{\poq{aq/b,aq/c}{n}}.\nonumber
\end{eqnarray}
Multiply both sides of this identity by
$\left(\frac{aq}{de}\right)^{n}\frac{\poq{d,e}{n}}{\poq{aq/d,aq/e}{n}}$
and reformulate the resulting identity to arrive at
\begin{eqnarray}&&\sum^{n}_{k=0}\frac{\poq {q^{-n},aq^{n},aq/bc,d,e}{k}} {\poq
{q,aq/b,aq/c,aq/d,aq/e}{k}}q^k\left(\frac{aq}{de}\right)^n
\frac{\poq{dq^{k},eq^{k}}{n-k}}{\poq{aq^{1+k}/d,aq^{1+k}/e}{n-k}}\nonumber\\
&&=\frac{\poq {b,c,d,e}{n}} {\poq {aq/b,aq/c,aq/d,aq/e}{n}}
\left(\frac{a^2q^2}{bcde}\right)^n.\label{identi}
\end{eqnarray}
Using two relations
\begin{eqnarray*}
&&\frac{\poq{q^{-k}}{i}}{\poq{deq^{-k}/a}{i}}
=\left(aq/de\right)^i\frac{\poq{q}{i}}{\poq{aq/(de)}{i}}\frac{\poq{q^{i+1}}{k-i}\poq{aq/(de)}{k-i}}{\poq{q}{k-i}\poq{aq^{i+1}/(de)}{k-i}}\quad
(k\geq i),
\\
&&\,_4\phi_3\left[%
\begin{array}{cccccccc}
  aq/(bc), & d, & e, & q^{-k} \\
  & aq/b, & aq/c,& deq^{-k}/a \\
\end{array};q,q
\right]\\
&&=\sum^{k}_{i=0}\frac{\poq {aq/bc,d,e}{i}} {\poq
{aq/(de),aq/b,aq/c}{i}}\left(\frac{aq}{de}\right)^i\frac{\poq{q^{i+1}}{k-i}\poq{aq/(de)}{k-i}}
{\poq{q}{k-i}\poq{aq^{i+1}/(de)}{k-i}},\end{eqnarray*} to reduces
(\ref{identi}), by exchanging the order of summations and
simplifying the resulting sums by the $q$-Pfaff-Saalsch\"{u}tz
formula, to
\begin{eqnarray}
\sum_{k=0}^n\frac{\poq{q^{-n}}{k}}{\poq{q}{k}}\frac{\poq{aq^{n}}{k}}{\poq{aq}{k}}q^k\,_4\phi_3\left[%
\begin{array}{cccccccc}
  aq/(bc), & d, & e, & q^{-k} \\
  & aq/b, & aq/c,& deq^{-k}/a \\
\end{array};q,q
\right] \nonumber\\
\times\frac{(aq,aq/(de);q)_{k}} {(aq/d,aq/e;q)_{k}}=\frac{\poq
{b,c,d,e}{n}} {\poq {aq/b,aq/c,aq/d,aq/e}{n}}
\left(\frac{a^2q^2}{bcde}\right)^n.\label{hyper}\end{eqnarray} Apply
the above
 Carlitz's  inversion  to (\ref{hyper}) and reformulate  the resulting identity in terms of standard hypergeometric
series as \begin{align}
&&\,_8\phi_7\left[%
\begin{array}{cccccccc}
  a, & q\sqrt{a}, & -q\sqrt{a}, & b, & c, & d,&e,&q^{-n} \\
   & \sqrt{a}, & -\sqrt{a}, & aq/b, & aq/c, & aq/d, & aq/e, & aq^{n+1}\\
\end{array};q,\frac{a^2q^{n+2}}{bcde}
\right]\nonumber \\&&=\frac{(aq,aq/(de);q)_{n}}
 {(aq/d,aq/e;q)_{n}}\,_4\phi_3\left[%
\begin{array}{cccccccc}
  aq/(bc), & d, & e, & q^{-n} \\
  & aq/b, & aq/c,& deq^{-n}/a \\
\end{array};q,q
\right]\label{5.14}
    \end{align}
 Now,  define
\begin{align}
   F(x)=\,_8\phi_7\left[%
\begin{array}{cccccccc}
  a, & q\sqrt{a}, & -q\sqrt{a}, & b, & c, & d,&e,&1/x \\
   & \sqrt{a}, & -\sqrt{a}, & aq/b, & aq/c, & aq/d, & aq/e, & aqx\\
\end{array};q,\frac{a^2q^{2}x}{bcde}
\right];\nonumber\end{align}
\begin{align} S(x)=\frac{(aq,aq/(de);q)_{\infty}}
{(aq/d,aq/e;q)_{\infty}}\frac{(aqx/e,aqx/d;q)_{\infty}}
{(aqx,aqx/(de);q)_{\infty}}\,_4\phi_3\left[%
\begin{array}{cccccccc}
  aq/(bc), & d, & e, & 1/x \\
  & aq/b, & aq/c, & de/(ax)  \\
\end{array};q,q
\right].\nonumber
    \end{align}
    It is immediately seen that $F(x)$ is analytic in open
set around $x=0$. So is $S(x)$ if we define particularly
$$S(0)=\frac{(aq,aq/(de);q)_{\infty}}
{(aq/d,aq/e;q)_{\infty}}\,_3\phi_2\left[%
\begin{array}{ccc}
  aq/(bc), & d, & e \\
   & aq/b, & aq/c \\
\end{array};q,\frac{a}{de}
\right]$$ and $ \poq{de/(ax)}{k}=0$ at finite points $x=deq^i/a$,
which can be guaranteed by the requirement that the $\,_4\phi_3$
series  be terminating. Further, Eq.(\ref{5.14}) means that for each
$n\geq 0$, $F(q^{n})=S(q^{n})$.
     Then by Theorem \ref{dl}, we have that
$ F(x)=S(x). $ Substitute $x$ by $1/f$  to get the desired result.
\qed

 For the limitation of space, we summarize in the following table some well-known summation and transformation formulas
 which can be shown by the same argument. Formula numbers in these formulas  refer
to the appendix of  Gasper and Rahman's book \cite{10}.
\begin{table}[!h]
\tabcolsep 0pt \caption{Table of $(f,g)-expansions$} \vspace*{-12pt}
\begin{center}
\def\temptablewidth{1\textwidth}
{\rule{\temptablewidth}{1pt}}
\begin{tabular*}{\temptablewidth}{@{\extracolsep{\fill}}c|c|c|c}
{\bf Sums or Transformations} & {\bf $(f,g)$-Expansions} &$x_n$ & $b_n$ \\
\hline
$q$-Gauss sum II.8 &$f=1,g=x-y$&$\setminus$&$q^n$\\
\hline
$q$-Kummer sum II.9  &$f=1-xy,g=x-y$&$aq^n$&$q^n$\\
\hline
 $q$-analogue of Bailey's $\,_2F_1$ sum II.10 &$f=g=x-y$&$aq^n$&$cq^n$\\
\hline
 $q$-Dixon sum $\,_4\phi_3$ II.13 &$f=1-xy,g=x-y$&$aq^n$&$q^n$\\
\hline
 II.18  &$f=(1-axy)(1-b\frac{x}{y})$&$cq^n$&$-q^n$\\
  &$g=(x-y)(1-\frac{b}{axy})$&&\\
\hline
 Bailey's formula $\,_3\varphi_3$ sum II.31 &&$q^n$&$q^n$\\
 Bilateral analogue of Dixon's sum II.32  &$f=1-xy,g=x-y$&$aq^n$&$q^n$\\
 $\,_3\varphi_3$ sum II.33  &&$aq^n$&$q^n$\\
\hline
Heine's transformation III.3 &$f=1,g=x-y$&$\setminus$&$cq^n$\\
\hline
Jackson's transformation III.4 &$f=1,g=x-y$&$\setminus$&$cq^n$\\
\hline
 III.9 &$f=1,g=x-y$&$\setminus$&$q^n$\\
\hline
 III.23 &$f=g=x-y$&$aq^n$&$q^n$\\
\hline
 Singh's transformation III.21  &$f=(1-xy)(1-qxy), g=x^2-y^2$&$cq^n$&$q^n$\\
\hline
 III.35&$f=g=x-y$&$aq^n$&$q^n$
\end{tabular*}
{\rule{\temptablewidth}{1pt}}
\end{center}
\end{table}
\subsection{From nonterminating summations to terminating summations}
The problem stated precisely is that given a nonterminating
summation of $q$-series, as reverse order of the ``Ismail's
argument", one can always expect a (perhaps new) terminating
summation applying the $(f,g)$-inversion. This method  can be
sketched as follows: consider that
\begin{equation*}{}_{r}\phi _{s}\left[\begin{matrix}a_{1},\dots ,a_{r}\\ b_{1},\dots ,b_{s}\end{matrix}
; q, z\right]\end{equation*} where $q$ is the base, $z$ is the
variable. Clearly, it contains  $r+s$ parameters. As it should be,
we might choose one of them as a new variable and the former
variable $z$ as a parameter such that it can be reformulated as the
form of (\ref{expandef3}). If succeed, then we are able to use the
$(f,g)$-expansion formula to derive a (new) terminating summation
formula. To illustrate, let take the $q$-binomial theorem as an
example.
\begin{lz} For $|q|<1$ and a variable $z: |z|<1$,
\begin{eqnarray}
  \sum_{n=0}^{\infty}\frac{(a;q)_n}{(q;q)_n}z^n=
  \frac{(az;q)_{\infty}}{(z;q)_{\infty}}\label{binom}
\end{eqnarray}
is equivalent to its finite form (\ref{binomfinite}), namely
\begin{eqnarray}
 \sum_{k=0}^{n}(-1)^kq^{\binom{k+1}{2}-nk}{n\brack
  k}_q\poq{z}{k}=z^n.\label{comform}
\end{eqnarray}
\end{lz}
{\sl Proof.} Now, we replace the parameter $a$ by a new variable $x$
and take $z$ as a parameter, and then reformulate (\ref{binom}) as
\begin{eqnarray*}
  \sum_{k=0}^{\infty}\frac{z^k}{(q;q)_k}(x;q)_k=
  \frac{(zx;q)_{\infty}}{(z;q)_{\infty}},
\end{eqnarray*}
which turns out to be equivalent to, under the substitution
$x\mapsto 1/x, z\mapsto zx$,
\begin{eqnarray}
  \sum_{k=0}^{\infty}\frac{(-1)^kz^k}{(q;q)_k}\prod_{i=0}^{k-1}(q^{i}-x)=
  \frac{(z;q)_{\infty}}{(zx;q)_{\infty}}.\label{binom00}
\end{eqnarray}
Letting $x=q^{n}$ yields a terminating  summation formula
\begin{eqnarray}
  \sum_{k=0}^{n}\frac{(-1)^kz^k}{(q;q)_k}\prod_{i=0}^{k-1}(q^{i}-q^{n})=
(z;q)_{n}\label{binom1}.
\end{eqnarray}
 Now, applying
 the $(1,x-y)$-inversion given in Theorem \ref{math4} to
(\ref{binom1}) to arrive at
\begin{eqnarray*}
 \frac{(-1)^nz^n}{(q;q)_n}=
 \sum_{k=0}^{n}\frac{(z;q)_{k}}{\displaystyle\prod_{i=0,i\neq k}^{n}(q^{i}-q^{k})}
 ,
\end{eqnarray*}
which reduces after simplification to the desired result.
Conversely, in the light of the $(1,x-y)$-inversion, (\ref{comform})
is equivalent to (\ref{binom1}), the latter states that
(\ref{binom00}) is valid for $x=q^n, n=0,1,2,\cdots$. Then by the
generalized ``Ismail's argument", (\ref{binom00}), i.e., the
$q$-binomial theorem holds.
 \qed

Clearly, (\ref{comform}) is just a $q$-analogue of the Newton
binomial formula
\begin{eqnarray}\sum_{k=0}^n\left(\!\!
\begin{array}{c}
  n \\
  k
\end{array}\!\!
\right)(1-y)^{k}y^{n-k}=1.\end{eqnarray}

 The next simple example is  how to apply our
argument  to transformations of $q$-series.

 \begin{lz}  Heine's transformation formula of ${}_{2}\phi _{1}$
series {\rm (cf.\cite[III.1]{10})}
\begin{eqnarray}
 && \,_2\phi_1\left[%
\begin{array}{cccccccc}
  a, & b \\
   & c  \\
\end{array};q,z
\right]=\frac{\poq{b}{\infty}(az;q)_{\infty}}{\poq{c}{\infty}(z;q)_{\infty}}{}_{2}\phi
_{1}\left[%
\begin{array}{cccccccc}
 z,& c/b\\
   & az
\end{array};q,b
\right] \label{heine1}\
\end{eqnarray}
is  equivalent to the following finite form of the $q$-binomial
theorem
\begin{eqnarray}
  \sum_{k=0}^{n}(-1)^kq^{\binom{k}{2}}{n\brack k}_qc^k=\poq{c}{n}.
\end{eqnarray}
\end{lz}
{\sl Proof.}
 In fact, take the parameter $b$ in (\ref{heine1}) as a
new variable $x$ and define
$$F(x)= \frac{\poq{c}{\infty}(z;q)_{\infty}}{\poq{x}{\infty}(az;q)_{\infty}}\,_2\phi_1\left[%
\begin{array}{cccccccc}
  a, & x \\
   & c  \\
\end{array};q,z
\right].$$ Thus, (\ref{heine1}) can be rewritten as
\begin{eqnarray}
 F(x)={}_{2}\phi
_{1}\left[%
\begin{array}{cccccccc}
 z,& c/x\\
   & az
\end{array};q,x
\right]. \label{heine2}\
\end{eqnarray}
So, it needs only to verify  that the r.h.s.\,of (\ref{heine2}) is
the $(1,x-y)$-expansion formula of $F(x)$  with parameters
$b_i=cq^{i}$ by calculating the $n$-th order $(1,x-y)$-difference of
$F(x)$. The result is
\begin{eqnarray*}
  \frac{\poq{aq^nz}{\infty}}{\poq{q^nz}{\infty}}c^n&=&\sum_{k=0}^{n}(-1)^kq^{\binom{k+1}{2}-nk}{n\brack k}_q\sum_{j\geq 0}
  \frac{\poq{a}{j}\poq{c}{j+k}}{\poq{q}{j}\poq{c}{j}}z^j.
\end{eqnarray*}
By employing the $q$-binomial theorem to expand
$\poq{aq^nz}{\infty}/\poq{q^nz}{\infty}$ in terms of $z^i$ and then
equating the coefficients of $z^m$ on both sides of this identity
leads to
\begin{eqnarray}
  \sum_{k=0}^{n}(-1)^kq^{\binom{k+1}{2}-nk}{n\brack k}_q\poq{cq^m}{k}=c^nq^{mn}.
\end{eqnarray}
Replace $cq^m$ by $c$. Then we get its simplified form:
\begin{eqnarray}
  \sum_{k=0}^{n}(-1)^kq^{\binom{k+1}{2}-nk}{n\brack
  k}_q\poq{c}{k}=c^n,\quad\mbox{i.e.,}\,\,\mathcal{D}_{q,x}^n\left\{\frac{\poq{c}{\infty}}{\poq{cx}{\infty}}\right\}|_{x=1}=c^n,
\end{eqnarray}
which in turn  is  equivalent to (\ref{binomfinite}), i.e., the
finite form of the $q$-binomial theorem.

\qed

 Performing the argument described in the above examples, we obtain corresponding terminating summations for some well-known
summation formulas of basic hypergeometric series
(cf.\cite[Appendixes I-III]{10}).
   \begin{lz}
    Heine's $q$-analogue of the Gauss summation
    formula {\rm (cf.\cite[II.8]{10})} (previously given as Example
\ref{lz})
    \begin{eqnarray}
     \,_2\phi_1\left[%
\begin{array}{cccccccc}
  a, & b \\
   & c  \\
\end{array};q,\frac{c}{ab}
\right]&=&
      \frac{(c/a,c/b;q)_{\infty}}{(c,c/(ab);q)_{\infty}}\label{gauss}
    \end{eqnarray}
    is equivalent to
\begin{equation}\sum_{k=0}^n(-1)^{k}{n\brack k}_qq^{\binom{n-k}{2}} \frac{\poq{c/a}{k}}{\poq
{c}{k}}=q^{\binom{n}{2}}\frac{\poq {a}{n}(c/a)^{n}}{\poq
{c}{n}}.\label{6.5}\end{equation}
\end{lz}
 As a byproduct, let $a\mapsto
\infty$ in (\ref{6.5}). It yields
\begin{equation}\sum_{k=0}^n(-1)^{n-k}{n\brack k}_qq^{\binom{n-k}{2}}\frac{1}{\poq
{c}{k}}=q^{2\binom{n}{2}}\frac{c^{n}}{\poq {c}{n}}.\end{equation}

 \begin{lz} Jackson's transformation formula of ${}_{2}\phi _{1}$
 to ${}_{2}\phi _{2}$ series {\rm (cf.\cite[III.4]{10})}
\begin{eqnarray}
 && \,_2\phi_1\left[%
\begin{array}{cccccccc}
  a, & b \\
   & c  \\
\end{array};q,z
\right]=\frac{(az;q)_{\infty}}{(z;q)_{\infty}}{}_{2}\phi
_{2}\left[%
\begin{array}{cccccccc}
 a,& c/b\\
  c, & az
\end{array};q,bz
\right]  \label{heine}\
\end{eqnarray}
is equivalent to
\begin{align} \sum_{k=0}^{n}(-1)^{n-k}{n\brack k}_qq^{\binom{k}{2}}\frac{1}{\poq{bq^{n-k}}{m+1}}={m+n\brack
m}_q \frac{b^nq^{\binom{n}{2}}\poq{q}{n}}{\poq{b}{m+n+1}},\quad
m\geq 0.\label{heineid}
\end{align}
\end{lz}
The details for the proofs of these examples are left to the
interested reader. As an interesting application of the expansion
formula (\ref{1.9}), now we can derive two new $q$-identities from
the following indefinite bibasic summation formula of Gasper
(cf.\cite[Eq.(1.14)]{8})
\begin{align}
\sum_{k=0}^{m}\frac{(1-ap^kq^{k})(1-bp^kq^{-k})}{(1-a)(1-b)}
\frac{(a,b;p)_k(x,a/(bx);q)_k}{(q,aq/b;q)_k(ap/x,bpx;p)_k}q^k\nonumber\\=\frac{(ap,bp;p)_m(xq,aq/(bx);q)_m}{(q,aq/b;q)_m
(ap/x,bpx;p)_m}.
\end{align}
by the same technique as above.
\begin{dl} For any integers $N, m\geq 0$, it holds
\begin{align}
\sum_{K=0}^N(-1)^{K}q^{\binom{K+1}{2}
}\begin{bmatrix}N\\K\end{bmatrix}_q
\frac{1-q^{m+1}}{1-q^{m+1+K}}\frac{1-aq^{2K+2m+2}/b}{1-aq^{m+1+K}/b}\nonumber
\\\frac{(a(pq)^{m+1}q^K,b(p/q)^{m+1}q^{-K};p)_{N}}
{(aq^{2m+2+K}/b;q)_{N+1}}=\frac{(ap^{m+1},bp^{m+1};p)_{N}}
{(aq^{m+1}/b;q)_{N+1}}/{N+m+1\brack m+1}_q.
\end{align}
In particular,
\begin{align}
\mathcal{D}_{q,x}^N\left\{\frac{(a(pq)^{m+1}x;p)_{N}}{1-xq^{m+1}}\right\}|_{x=1}=\frac{q^{N(m+1)}(ap^{m+1};p)_{N}}{1-q^{m+1}}/{N+m+1\brack
m+1}_q.
\end{align}
\end{dl}
{\sl Proof.} Observe that Gasper's indefinite summation formula can
be restated shortly as $F(x)=S(x)$ by choosing
$f(x,y)=(1-axy)(1-b\frac{x}{y}),g(x,y)=(x-y)(1-\frac{b}{axy})$, in
this case, $f(x,y)\in \mbox{\sl Ker}\mathcal{L}^{(g)}_{3}$, and
 $b_i=q^i, x_i=p^i$, $|q|<1$, as well as by defining
\begin{eqnarray}
&&F(x)=\frac{(ap,bp;p)_m(q/x,aqx/b;q)_m}{(q,aq/b;q)_m(apx,bp/x;p)_m};\\
&& S(x)=\sum_{k=0}^{\infty}G(k)f(x_k,b_k)
\frac{\prod_{i=0}^{k-1}g(b_i,x)}{\prod_{i=1}^{k}f(x_i,x)},\label{check1}
\end{eqnarray}
where the coefficients
\begin{eqnarray*}
G(n)=\left\{%
\begin{array}{ll}
    0, & \hbox{for}\,\, n\geq m+1; \\
   \displaystyle\frac{(a,b;p)_na^nq^{n(n+1)/2}}{(1-a)(1-b)(q,aq/b;q)_nb^n}, &\hbox{for}\,\, n\leq m.
 \end{array}%
\right.
   \end{eqnarray*}
In the sense of the representation problem of $(f,g)$-series, $S(x)$
is just the $(f,g)$-expansion formula of $F(x)$, since $F(x)$ is
rational in $x$, thus analytic around $x=0$, so is $S(x)$ for the
series in the r.h.s.\,of (\ref{check1}) is a finite sum. Thus, by
the uniqueness of the $(f,g)$-expansion formula, namely, Lemma
\ref{uni}, we have that for $n\geq m+1$ (the result corresponding to
the case $n\leq m$ is trivial),
\begin{eqnarray*}
&&
G(n)=\,{\bf\mathbb{D}}^{(n)}_{\left((1-axy)(1-b\frac{x}{y}),(x-y)(1-\frac{b}{axy})\right)}\left[\begin{matrix}1,q,q^2,
\dots,q^{n}\\p,p^2,\dots,p^{n-1}\end{matrix}\right]\{F(x)\}
\nonumber\\
&=&\frac{1}{(q;q)_n}
\sum_{k=0}^n(-1)^{n-k}q^{\binom{k+1}{2}-nk}{n\brack
k}_q(1-\frac{b}{a}q^{-2k})\frac{(apq^k,bpq^{-k};p)_{n-1}}
{(\frac{b}{a}q^{-k};q^{-1})_{n+1}}F(q^k)=0.
   \label{check}
\end{eqnarray*}
Putting the values of $F(q^k)$ into it and simplifying the resulted,
we obtain that
\begin{align}
\sum_{k=m+1}^n(-1)^{k-m-1}q^{\binom{k}{2}+\binom{m+1}{2}-mk}\begin{bmatrix}n\\k\end{bmatrix}_q\begin{bmatrix}k-1\\m\end{bmatrix}_q
\frac{1-aq^{2k}/b}{1-aq^k/b}\nonumber
\\\frac{(ap^{m+1}q^k,bp^{m+1}q^{-k};p)_{n-m-1}}
{(aq^{m+k+1}/b;q)_{n-m}}=\frac{(ap^{m+1},bp^{m+1};p)_{n-m-1}}
{(aq^{m+1}/b;q)_{n-m}}.\label{3.67}
\end{align}
Next, apply the substitution  $n\mapsto N+m+1$ and $k\mapsto K+m+1$
to (\ref{3.67}), and then use the basic relation
$$
\begin{bmatrix}n\\k\end{bmatrix}_q\begin{bmatrix}k-1\\m\end{bmatrix}_q={n\brack
m+1}_q{n-m-1\brack k-m-1}_q\frac{1-q^{m+1}}{1-q^{k}}
$$
to simplify the resulting identity. The final result is
\begin{align}
{N+m+1\brack
m+1}_q\sum_{K=0}^N(-1)^{K}q^{\binom{K+1}{2}}\begin{bmatrix}N\\K\end{bmatrix}_q
\frac{1-q^{m+1}}{1-q^{m+1+K}}\frac{1-aq^{2K+2m+2}/b}{1-aq^{K+m+1}/b}\nonumber
\\\frac{(a(pq)^{m+1}q^K,b(p/q)^{m+1}q^{-K};p)_{N}}
{(aq^{2m+2+K}/b;q)_{N+1}}=\frac{(ap^{m+1},bp^{m+1};p)_{N}}
{(aq^{m+1}/b;q)_{N+1}}.\label{3.68}
\end{align}
If we divide both sides of (\ref{3.68}) by $b^{N+1}$ and then take
the limit $b\mapsto 0$ on both sides, then we obtain further
\begin{align}{N+m+1\brack
m+1}_qq^{-N(m+1)}\sum_{K=0}^N(-1)^{K}q^{\binom{K+1}{2}-NK}\begin{bmatrix}N\\K\end{bmatrix}_q\nonumber\\
\frac{1-q^{m+1}}{1-q^{m+1+K}}(a(pq)^{m+1}q^K;p)_{N}=
(ap^{m+1};p)_{N}.
\end{align}
As stated in Proposition \ref{yl11}, it can be reformulated in term
of $\mathcal{D}_{q,x}$ as the desired form. \qed
\section{Conclusions}
We hope that  the generalized ``Ismail's argument" or the
representation of analytic functions in terms of $(f,g)$-series in
this article
 is  a new general  approach to the basic  hypergeometric series.
Perhaps, the most intriguing case is
  that  the $n$th order $(f,g)$-difference  of
$F(x)\in {\mathcal H}(\Omega)$ $$\,{\bf\mathbb{D}}^{(n)}_{(f,g)}\left[\begin{matrix}b_0,b_{1},\dots ,b_{n}\\
x_{1},\dots ,x_{n-1}\end{matrix}\right]\{F(x)\}$$ can be evaluated
in a closed form. If so, it provides an affirmative  answer to a
problem poised by Marco and  Parcet in \cite[\S 5]{mar}. On the
other hand, the $(f,g)$-expansion
 formula of $F(x)$ is in fact a rational approximation to $F(x)$ if $f(x,y)$ and
$g(x,y)$ are polynomials of two variables  $x$ and $y$. Thus, it is
necessary
 to study application of  this expansion formula as well as the $n$-th order $(f,g)$-difference operator in the theory of (numerical) approximation.
  The
same problem remains open to the  expansion formula (\ref{1.10}) in
the theory of elliptic hypergeometric series \cite{1000}. Besides,
we believe our results are also partial solutions to a
 problem of Gessel and Stanton proposed in \S 9 of \cite{111}.
 All these problems will be discussed
 in our forthcoming papers.
\section{Appendix}
In this appendix, the proofs of the crucial fact
 that $\lim_{n\mapsto \infty}S(n,x)/a_n=0$ invoked in Theorem
\ref{dl} and the fact that Gessel and Stanton's $q$-analogue
(\ref{gessel1})/(\ref{gessel2}) is equivalent to the
$(1-xy,x-y)$-expansion formula with respect to geometric series are
displayed in details.

\begin{yl}\label{bu2}Preserve all assumptions as in Theorem \ref{dl} and let  $K_{n,k}(x)$ be given  by (\ref{knk}). Define that
\begin{eqnarray*}
S(n,x)=\sum_{k=1}^{n-1}(-1)^{n-k}A^{k}p^{\binom{k+1}{2}} {n-1\brack
k}_p\frac{\poq{q}{n-k}}{\poq{q}{n}}
\sum_{i=0}^k(-1)^iq^{(n-k+1)i+\binom{i}{2}}{k\brack
i}_qK_{n-k,n-k}(q^ix).
\end{eqnarray*}
Then \begin{eqnarray} \lim_{n\mapsto \infty}
S(n,x)/a_n=0\end{eqnarray}  for  $x=1,q$.
\end{yl}
\pf We proceed to  show the desired result by Tannery's theorem,
more precisely, to find  $C(n)$ and $T_k$ such that
$|t_{n,k}/C(n)|<T_k,$ $\sum_{k=0}^\infty T_k$ is convergent, and for
each fixed $k$, there holds $\lim_{n\mapsto\infty}t_{n,k}/C(n)$.
Obviously, $S(n,x)$ can be reformulated as
\begin{eqnarray*}
&&S(n,x)=\sum_{k=1}^{n-1}t_{n,k};\\
&&t_{n,k}=(-1)^{k}A^{n-k}p^{\binom{n-k+1}{2}} {n-1\brack
k-1}_p\frac{\poq{q}{k}}{\poq{q}{n}}\sum_{i=0}^{n-k}(-1)^iq^{(k+1)i+\binom{i}{2}}{n-k\brack
i}_qK_{k,k}(q^{i}x).
\end{eqnarray*}
These required functions can be found by virtue of the following
inequalities.

\textbf{(a) }Using the basic relation of the  $q$-binomial
coefficients
$$
{n+1\brack i}_q=q^i{n\brack i}_q+{n\brack i-1}_q
$$
and induction on $i\leq n$, we obtain at first that
\begin{eqnarray}
\left|{n\brack i}_q \right|\leq {n\brack i}_{|q|}.
\end{eqnarray}
On the other hand,  from (\ref{binomfinite}) it follows
\begin{eqnarray*}
\sum_{i=0}^{n-k}q^{\binom{i}{2}}{n-k\brack i}_qx^i=\poq{-x}{n-k}.
\end{eqnarray*}
Combined these two relations. It is easily seen that
\begin{eqnarray}
\sum_{i=0}^{n-k}|q|^{(k+1)i+\binom{i}{2}}\left|{n-k\brack i}_q
\right|\leq \sum_{i=0}^{n-k}|q|^{(k+1)i+\binom{i}{2}}{n-k\brack
i}_{|q|}=\frac{\poqt{-|q|}{n}}{\poqt{-|q|}{k}}. \label{ineq1}
\end{eqnarray}

\textbf{(b)} Note that the real-valued function of a variable $x$
\begin{eqnarray*}
\bar{K}_{k,k}(x)=\sum_{r=0}^{\infty}\left|a_{r+k}\right|{r+k\brack
r}_{|q|}|x|^{r}
\end{eqnarray*}
 is nondecreasing for $|x|\leq \min\{1,R\}$, which gives that
\begin{eqnarray}\bar{K}_{k,k}(q^{i}x)\leq \bar{K}_{k,k}(x)\leq \bar{K}_{k,k}(1).\label{ineq2}
\end{eqnarray}

\textbf{(c}) Given $|p|<1$, the real-valued function of degree two
in a variable $x$ ($t$  is as a parameter)
$$
y(x)=(x-k)(x-k+1)/2\, \ln |p|+(x-k)\ln |A|- t\ln |c_0| x
$$
on the interval $x\geq k$ takes its the maximum value $M_k$ at the
point $x=k-1/2+(t\ln |c_0|-\ln |A|)/\ln |p|$, under the subsidiary
condition that $\ln m \leq t\leq \ln M$,  and
\begin{eqnarray*}
&&M_k=(c-1/2) \ln(|A/c_0|)+(c^2/2-1/8)\ln |p|-t\, k \ln |c_0|
\end{eqnarray*}
where $c=(t\ln |c_0|-\ln |A|)/\ln |p|,$ since the second derivative
$y^{''}(x)=\ln |p|<0$. It ensures that there must exist a constant
$\bar{M}$ independent of $k$ such that
\begin{eqnarray}
|A^{n-k}p^{\binom{n-k+1}{2}}/a_n|=\exp(y(n))\leq
\bar{M}R^k.\label{ineq3}
\end{eqnarray}
 So, with the help of  Inequalities (\ref{ineq1})-(\ref{ineq3}), we  obtain that
\begin{eqnarray}
|t_{n,k}|&=&\left|(-1)^{k}A^{n-k}p^{\binom{n-k+1}{2}} {n-1\brack
k-1}_p\frac{\poq{q}{k}}{\poq{q}{n}}\sum_{i=0}^{n-k}(-1)^iq^{(k+1)i+\binom{i}{2}}{n-k\brack
i}_qK_{k,k}(q^{i}x)\right|\nonumber\\
 &&\leq
\bar{M}R^k|a_n|\times {n-1\brack
k-1}_{|p|}\times\frac{|\poq{q}{k}|}{|\poq{q}{n}|}\times\bar{K}_{k,k}(x)\times
\sum_{i=0}^{n-k}|q|^{(k+1)i+\binom{i}{2}}\left|{n-k\brack i}_q
\right|\nonumber\\
&&\leq \bar{M}R^k|a_n|\times
\frac{\poqpt{|p|}{n-1}|\poq{q}{k}|}{\poqpt{|p|}{k-1}\poqpt{|p|}{\infty}|\poq{q}{n}|}\times\bar{K}_{k,k}(1)\times
\frac{\poqt{-|q|}{n}}{\poqt{-|q|}{k}}=|C(n)|T_k,
\end{eqnarray}
where, for simplicity, we define
\begin{eqnarray*}
&&C(n)=\frac{\bar{M}
a_n}{\poqpt{|p|}{\infty}}\frac{\poqpt{|p|}{n-1}\poqt{-|q|}{n}}{|\poq{q}{n}|};\\
&&T_k=A_kB_k,
A_k=\frac{|\poq{q}{k}||a_k|R^k}{\poqpt{|p|}{k-1}\poqt{-|q|}{k}},B_k=\frac{\bar{K}_{k,k}(1)}{|a_k|}.
\end{eqnarray*} Note that in the above
estimate, we utilize an inequality
$$
\poqpt{|p|}{n-k}\geq \poqpt{|p|}{\infty}\qquad\mbox{for}\quad |p|<1.
$$ Under the known conditions, we see that
$\lim_{k\mapsto\infty} B_k=1/|(1/R;|q|)_{\infty}|$ and $\sum_{k\geq
0}A_k$ is convergent. Now, by Cauchy criterion for convergence, it
is easily found that $\sum_{k\geq 0}T_k$ is convergent. Now, we are
in a position to apply
 Tannery's theorem to the sum $S(n,x)/C(n)$. The result is
\begin{eqnarray*}
\lim_{n\mapsto\infty}\frac{S(n,x)}{C(n)}&=&\frac{|\poq{q}{\infty}|}
{\bar{M}\poqt{-|q|}{\infty}}\sum_{k=1}^{\infty}\lim_{n\mapsto\infty}(-1)^{k}\frac{A^{n-k}p^{\binom{n-k+1}{2}}}{a_n}
{n-1\brack
k-1}_p\frac{\poq{q}{k}}{\poq{q}{n}}\\
&\times&\lim_{n\mapsto\infty}\sum_{i=0}^{n-k}(-1)^iq^{(k+1)i+\binom{i}{2}}{n-k\brack
i}_qK_{k,k}(q^{i}x)=0
\end{eqnarray*}
 by using the basic relations
\begin{eqnarray*}
&&\lim_{n\mapsto\infty}\frac{(-1)^{k}A^{n-k}p^{\binom{n-k+1}{2}}}{a_n}
=0,\,\, \lim_{n\mapsto\infty}\frac{K_{n,n}(x)}{
a_n}=\frac{1}{\poq{xc_0}{\infty}};\\
&&\lim_{n\mapsto\infty}\sum_{i=0}^{n-k}(-1)^iq^{(k+1)i+\binom{i}{2}}{n-k\brack
i}_qK_{k,k}(q^{i}x)=\frac{\poq{q}{\infty}}{\poq{q}{k}}\sum_{r=0}^{\infty}\frac{a_{r+k}}{\poq{q}{r}}x^r.
\end{eqnarray*}
The last limitation is also a consequence of applying Tannery's
theorem. Finally, a simple observation that
\begin{eqnarray*}
\lim_{n\mapsto\infty}\frac{S(n,x)}{C(n)}=\frac{|\poq{q}{\infty}|}
{\bar{M}\poqt{-|q|}{\infty}}\lim_{n\mapsto\infty}\frac{S(n,x)}{a_n}
\end{eqnarray*}
 gives the complete proof of the claimed.\qed

As a  previously unknown fact, as mentioned earlier, that the
$(1-xy,x-y)-$expansion formula with respect to the parameter
sequences $b_i=q^i, x_i=Ap^i$ is actually equivalent to  Gessel and
Stanton's $q$-analogue (\ref{gessel1})/(\ref{gessel2}) of the
Lagrange inversion formula for $y=x/(1-x)^{b+1} $(cf.\cite[p.180,
Theorem 3.7]{111}). Now it can be made clear by reducing the
following matrix version of their $q$-analogue from the
$(1-xy,x-y)-$expansion formula in Theorem \ref{dl}.
\begin{yl}\label{yl112}
 Let $B=(B_{n,k})$ and $B^{-1}=(B^{-1}_{n,k})$ be inverses of each other, where
$$
B_{n,k}=\frac{(Ap^kq^k;p)_{n-k}}{\poq{q}{n-k}}q^{-nk}.
$$
Then
$$
B_{n,k}^{-1}=(-1)^{n-k}q^{\binom{n-k+1}{2}+nk}
\frac{(1-Ap^kq^k)(Aq^np^{n-1};p^{-1})_{n-k-1}}{\poq{q}{n-k}}.
$$
\end{yl}
{\sl Proof.} It suffices to, by Theorem \ref{dl},  calculate the
$n$th order $(1-xy,x-y)$-difference of $F(x)$
\begin{eqnarray}
&&G(n)=\sum_{k=0}^{n} F(q^k)\frac {\prod_{i=1}^{k-1} (1-Ap^iq^k)}
 {\prod_{i=0}^{k-1}(q^i-q^k)}
 \frac{\prod_{i=k}^{n-1}(1-Ap^iq^{k})}
 {\prod_{i=k+1}^{n}(q^i-q^k)}\label{id3}.
\end{eqnarray}
Note that
\begin{align*}
 \prod_{i=k}^{n-1}(1-Ap^iq^{k})=\poq{Ap^kq^{k}}{n-k},\,\,\prod_{i=k+1}^{n}(q^i-q^k)=(-1)^{n-k}q^{k(n-k)}\poq{q}{n-k}.
\end{align*}
By these notes, (\ref{id3}) can be reformulated as
\begin{align}
(-1)^nG(n)=\sum_{k=0}^{n}B_{n,k} \left\{ (-1)^kq^{\binom{k+1}{2}}
\frac{(Apq^{k};p)_{k-1}}
 {\poq{q}{k}}F(q^k)\label{id4}\right\}.
\end{align}
On the other hand, taking $x=q^{n}$ in the $(1-xy,x-y)$-expansion
formula yields
\begin{eqnarray*}
F(q^n)&=&\sum_{k=0}^{n} (1-Ap^kq^{k})
 \frac{\prod_{i=0}^{k-1}(q^i-q^n)}
{\prod_{i=1}^{k}(1-Ap^iq^{n})}G(k)\\
&=&\frac{\prod_{i=0}^{n-1}(q^i-q^n)}{\prod_{i=1}^{n-1}(1-Ap^iq^{n})}\sum_{k=0}^{n}
(1-Ap^kq^{k}) \frac{\prod_{i=k+1}^{n-1}(1-Ap^iq^{n})}
{\prod_{i=k}^{n-1}(q^i-q^n)}G(k).
\end{eqnarray*}
A slight simplification gives
\begin{eqnarray*}
\frac{\prod_{i=1}^{n-1}(1-Ap^iq^{n})}{\prod_{i=0}^{n-1}(q^i-q^n)}F(q^n)
&=&\sum_{k=0}^{n} (1-Ap^kq^{k})
\frac{\prod_{i=k+1}^{n-1}(1-Ap^iq^{n})}
{\prod_{i=k}^{n-1}(q^i-q^n)}G(k),
\end{eqnarray*}
which becomes after rewritten in terms of $q$-shifted factorials
\begin{eqnarray*}
q^{-\binom{n}{2}}\frac{(Apq^{n};p)_{n-1}}{\poq{q}{n}}F(q^n)
&=&\sum_{k=0}^{n} (1-Ap^kq^{k}) q^{-\binom{n}{2}+\binom{k}{2}}
\frac{(Aq^np^{n-1};p^{-1})_{n-k-1}} {\poq{q}{n-k}}G(k).
\end{eqnarray*}
Finally, we get
\begin{eqnarray}
&&\label{add1}\left\{(-1)^nq^{\binom{n+1}{2}}\frac{(Apq^{n};p)_{n-1}}{\poq{q}{n}}F(q^n)\right\}\\
&=&\sum_{k=0}^{n}(-1)^{n-k} (1-Ap^kq^{k}) q^{\binom{n-k+1}{2}+ nk}
\frac{(Aq^np^{n-1};p^{-1})_{n-k-1}}
{\poq{q}{n-k}}\left\{(-1)^kG(k)\right\}.\nonumber
\end{eqnarray}
 Define that
$$
\left\{%
\begin{array}{ll}
    &f_n=(-1)^nG(n); \\
&\\
    &a_n=(-1)^nq^{\binom{n+1}{2}}
    \frac{(Apq^{n};p)_{n-1}}
 {\poq{q}{n}}F(q^n),
\end{array}%
\right.
$$
and $X=(f_0,f_1,f_2,\cdots,f_n,\cdots)^T,
Y=(a_0,a_1,a_2,\cdots,a_n,\cdots)^T$, the superscript $T$ denotes
the transpose of  matrix. With these notation, (\ref{id4}) and
(\ref{add1}) can be reformulated respectively as
\begin{eqnarray*}
   f_n&=&\sum_{k=0}^nB_{n,k}a_k\Leftrightarrow X=BY;\\
  a_n&=&\sum_{k=0}^{n}(-1)^{n-k} (1-Ap^kq^{k}) q^{\binom{n-k+1}{2}+
nk} \frac{(Aq^np^{n-1};p^{-1})_{n-k-1}}
{\poq{q}{n-k}}f_k\Leftrightarrow Y=B^{-1}X.
\end{eqnarray*}
From this one can read off that
$$
B_{n,k}^{-1}=(-1)^{n-k} (1-Ap^kq^{k}) q^{\binom{n-k+1}{2}+ nk}
\frac{(Aq^np^{n-1};p^{-1})_{n-k-1}} {\poq{q}{n-k}}.
$$
This gives the complete proof of the desired.\qed
\bibliographystyle{amsplain}

\begin{thebibliography}{99}

\bibitem{ahl} L. V. Ahlfors, Complex Analysis, 3rd ed., McGraw-Hill,
1978.
\bibitem{andrew} G. E. Andrews, {\em Identities in Combinatorics. II:
A q-analogue of the Lagrange inversion theorem}, Proc.Amer. Math.
Soc. \textbf{53}(1975), 240-245.


\bibitem{andrews}
\bysame, R. Askey, and R. Roy, Special Functions, Encyclopedia of
Mathematics and Its Applications, Vol. \textbf{71}, Cambridge
University Press, Cambridge, UK, 1999.

\bibitem{ask}
R. Askey and M.E.H. Ismail,\,{\sl A very-well-poised
$\,_6\varphi_6$}, Proc. Amer. Math. Soc.,\textbf{77}(1979), 218-222.

\bibitem{macd}G. Bhatnagar, {\sl Generalized bibasic series, and their $U(n)$ extensions},
Advances in Math. \textbf{131} (1997), 188-252.

\bibitem{boa} R. P. Boas, \,{\sl Tannery's theorem}, Math. Mag., 1965, \textbf{38}: 66.

\bibitem{4}
D. M. Bressoud, {\em Some identities for terminating q--series},
Math. Proc. Cambridge Philos. Soc. \textbf{89} (1981), 211--223.

\bibitem{5}
\bysame, {\em A matrix inverse}, Proc. Amer. Math. Soc. \textbf{88}
(1983), 446--448.

\bibitem{tttsss}
R. L. Burden, J. D. Faires, Numerical Analysis,  7th Edition,
Pacific Grove, CA: Brooks/Cole, 2001.

\bibitem{carliz1} L. Carlitz, {\em Some $q$-identities related to the theta
functions,} Bullettino della Unione Matematica Italiana
\textbf{17}(3) (1962), 172-178.

\bibitem{6}
\bysame, {\em Some inverse relations}, Duke Math. J. \textbf{40}
(1973), 893--901.

\bibitem{carliz}
\bysame, {\em Some q-expansion theorems}, Glas. Math. Ser. III
\textbf{ 8}(28) (1973), 205-214.

\bibitem{chu}W.Chu, L. Di Claudio, Classical Partition Identities and Basic Hypergeometric Series,
Quaderni del Dipartimento di Matematica dell'Universit¨¤ di Lecce;
\textbf{6} (2004), 149-150.


\bibitem{fine}N. J. Fine, Basic Hypergeometric Series and Applications,
in Mathematical Surveys and Monographs, Vol.\textbf{27}, Amer. Math.
Soc., Providence, 1988.


\bibitem{fu1}
Amy M. Fu and Alain Lascoux, {\em q-Identities from Lagrange and
Newton interpolation}, Advances in Appl. Math. \textbf{31} (2003),
527-531.

\bibitem{fu2}
\bysame, {\em Rational interpolation and basic hypergeometric
series}, arXiv:math.CO/0404063.

\bibitem{gar1} A. Garsia, {\em A $q$-analogue of the Lagrange inversion formula}, Houston J. Math.
\textbf{7} (1981), 205-237.

\bibitem{gar2}A. Garsia and J. Remmel,{\em A novel form of $q$-Lagrange inversion}, Houston J. Math. \textbf{12}
(1986), 503-524.

\bibitem{8}
G. Gasper, {\em Summation, transformation, and expansion formulas
for bibasic series}, Trans. Amer. Math. Soc. \textbf{312} (1989),
257--278.

\bibitem{188}
\bysame, {\em Elementary derivations of summation and transformation
formulas for $q$-series}, in Special Functions, $q$-Series and
Related Topics (M.E.H.Ismail, D.R.Masson and M.Rahman, eds),
Amer.Math.Soc., Providence, R.I., Fields Institute Communications
\textbf{14} (1997), 55-70.

\bibitem{rahman}
\bysame and M. Rahman, {\em An indefinite bibasic summation formula
and some quadratic, cubic and quartic summation and transformation
formulas}, Canad. J. Math. \textbf{42} (1990), 1-27.


\bibitem{10}
\bysame and M. Rahman, Basic Hypergeometric Series, second edition,
Cambridge University Press, Cambridge, 2004.


\bibitem{111}
I. Gessel and D. Stanton, {\em Application of q-Lagrange inversion
to basic hypergeometric series}, Trans. Amer. Math. Soc. \textbf{
277} (1983), 173-203.

\bibitem{16}
H. W. Gould and L. C. Hsu, {\em Some new inverse series relations},
Duke Math. J. \textbf{40} (1973), 885-891.

\bibitem{ismail0} M.E.H.Ismail, {\em A simple proof of Ramanujan's ${}_{1}\varphi_{1}$ sum,
} Proc. Amer. Math. Soc. \textbf{63} (1977), 185-186.


\bibitem{ismail1} \bysame and D.Stanton, {\em Applications of q-Taylor
theorems}, J. Comp. Appl. Math. \textbf{153} (2003), 259-272.

\bibitem{ismail2} \bysame, {\em $q$-Taylor theorems, polynomial expansions,
and interpolation of entire functions}, J. Approx. Theory
\textbf{123} (2003), 125-146

\bibitem{jackson}F.H. Jackson, {\em On q-functions and a certain difference
operator,} Trans. Roy Soc. Edin. \textbf{46} (1908), 253-281.

\bibitem{kratt} C. Krattenthaler, {\em A new q-Lagrange formula and some
applications}, Proc. Amer. Math. Soc. \textbf{90} (1984), 338--344.

\bibitem{18}
\bysame, {\em Operator methods and Lagrange inversion, a unified
approach to Lagrange formulas}, Trans. Amer. Math. Soc. \textbf{305}
(1988), 431--465.

\bibitem{99}\bysame,\ \textit{A new matrix inverse},\,Proc.Amer.Math.Soc.\textbf{124}\
(1996), 47-59.
\bibitem{koo}Tom H. Koornwinder, {\em Some simple applications and variants of the $q$-binomial
formula}, Informal note, 1999.

\bibitem{ttt}A. Lascoux and M. P. Sch\"{u}tzenberger, {\em  Symmetrization operators on
polynomial rings}, Funkt. Anal. \textbf{21} (1987), 77-78 .

\bibitem{liu}
Z. G. Liu, {\em An expansion formula for q-series and applications},
The Ramanujan J., \textbf{6}(2002), 429-447.

\bibitem{mar} B. L\'{o}pez, J. M. Marco and J. Parcet,
{\em Taylor series for the Askey-Wilson operator and classical
summation formulas}, Proc. Amer. Math. Soc. \textbf{134} (2006),
2259-2270.

\bibitem{0020}X. Ma,\ \textit{An extension of Warnaar's matrix inversion},\
Proc. Amer. Math. Soc., \textbf{133} (2005), 3179-3189.

\bibitem{0022}\bysame,\ \textit{Two finite forms of Watson's quintuple product identity and matrix
inversion}, Electron. J. Comb. \textbf{13}(1) (2006)\, \#R52, 8pp.

\bibitem{0021}\bysame,\ \textit{The $(f,g)$-inversion formula and
its applications: the $(f,g)$-summation formula}, Advances in Appl.
Math.,\textbf{37} (2006), to appear.

\bibitem{milne} S. C. Milne and G. Bhatnagar,\ \textit{A characterization of inverse relations}, Discrete
Math. \ \textbf{193}\ (1998),\ 235--245.

\bibitem{19}
M. Rahman, {\em Some quadratic and cubic summation formulas for
basic hypergeometric series}, Canad. J. Math. \textbf{45} (1993),
394--411.


\bibitem{20}
\bysame,  {\em Some cubic summation formulas for basic
hypergeometric series}, Utilitas Math. \textbf{36} (1989), 161--172.


\bibitem{3rr}
S. M. Roman, {\em  More on the umbral calculus, with emphasis on the
q-umbral calculus}, J. Math. Anal. Appl. 107 (1985), 222-254.

\bibitem{ryd} F. Ryde,  A contribution to the theory of linear
homogeneous geometric difference equations ($q$-difference
equations), Dissertation, Lund, 1921.

\bibitem{schlosser1}  M. Schlosser,\ \textit{A simple proof of Bailey's very-well-poised
$\,_6\psi_6$ summation}, Proc. Amer. Math. Soc., \textbf{130}
(2002), 1113-1123.

\bibitem{schlosser2}
\bysame, {\sl Inversion of bilateral basic hypergeometric series},
Electron. J. Comb. \textbf{10} (2003)\, \#R10, 27pp.

\bibitem{22}
D. Singer, {\em $q$-analogues of Lagrange inversion}, Ph.D. Thesis,
Univ. of California, San Diego, CA, 1992.

\bibitem{29} V.P.Spiridonov,\ \textit{Theta hypergeometric series}, in V. A.
Malyshev and A. M. Vershik (eds.), Asymptotic Combinatorics with
Applications to Mathematical Physics, 307-327, Kluwer Acad.Publ.,
Dordrecht, 2002.

\bibitem{11}
 D. Stanton, {\em Recent results for the q--Lagrange
inversion formula}, Ramanujan Revisited, ed. by Askey, Berndt,
Ramanathan, Rankin, Academic Press, 1988, 525-536.

\bibitem{ste} J. F. Steffensen, {\sl On divided differences,} Danske Vid. Selsk.
Math.-Fys. Medd., \textbf{17} (1939), 1-12.


\bibitem{1000} S. O. Warnaar,\ {\sl Summmation and transformation formulas for elliptic hypergeometric series,}
Constr. Approx.\textbf{18} (2002), 479--502.

\bibitem{1001}E. T. Whittaker, G. N. Watson,\ A Course of Modern Analysis,
Reprint of the fourth (1927) edition, Cambridge University Press,
Cambridge, 1996. \MR{97k:}01072

\bibitem{jiang}
J. Zeng, {\em On some q-identities related to divisor functions},
 Advances in Appl.
Math. \textbf{34} (2005), 313-315.


\end{thebibliography}

\end{document}